\documentstyle[11pt,bezier,epsfig]{article}
\newcommand{\nc}{\newcommand}
\nc{\rnc}{\renewcommand}
\rnc{\thesection}{\arabic{section}\setcounter{equation}{0}}
\rnc{\theequation}{\arabic{section}.\arabic{equation}}
\nc{\be}{\begin{equation}}
\nc{\ee}{\end{equation}}
\nc{\bea}{\begin{eqnarray}}
\nc{\eea}{\end{eqnarray}}
\nc{\bean}{\begin{eqnarray*}}
\nc{\eean}{\end{eqnarray*}}
\nc{\nn}{\nonumber}
\nc{\mb}{\mbox}
\nc{\ch}{\mb{ch}}
\nc{\sh}{\mb{sh}}
\nc{\tnh}{\mb{th}}
\nc{\Y}{{\cal Y}}
\nc{\T}{{\cal T}}
\nc{\st}{{\cal S}}
\newtheorem{th}{Theorem}
\newtheorem{lem}{Lemma}

\newtheorem{cj}{Conjecture}
\def\vertex[#1,#2,#3,#4,#5]{
 \begin{picture}(22,10)(0,10)
 \put(4,10){\line(1,0){12}}
 \put(10,4){\line(0,1){12}}
 \put(10,10){\oval(2.5,2.5)[bl]}
 \put(0,9){{\scriptsize $#1$}}
 \put(9,1.5){{\scriptsize $#2$}}
 \put(16.5,9){{\scriptsize $#3$}}
 \put(9,16.9){{\scriptsize $#4$}}
 \put(7,7){{\scriptsize $#5$}}
 \end{picture}}
\def\BWa{
 \begin{picture}(22,10)(0,10)
 \multiput(4,10)(1.8,0){7}{\line(1,0){1.1}}
 \put(10,4){\line(0,1){12}}
 \put(10,10){\oval(2.5,2.5)[bl]}
 \put(2,9){{\scriptsize $j$}}
 \put(9,1.5){{\scriptsize $\pm 1$}}
 \put(16.5,9){{\scriptsize $j$}}
 \put(9,16.9){{\scriptsize $\pm 1$}}
 \put(7,7){{\scriptsize $v$}}
 \end{picture}}
\def\BWb{
 \begin{picture}(22,10)(0,10)
 \multiput(4,10)(1.8,0){7}{\line(1,0){1.1}}
 \put(10,4){\line(0,1){12}}
 \put(10,10){\oval(2.5,2.5)[br]}
 \put(2,9){{\scriptsize $j$}}
 \put(9,1.5){{\scriptsize $\pm 1$}}
 \put(16.5,9){{\scriptsize $j$}}
 \put(9,16.9){{\scriptsize $\pm 1$}}
 \put(11,7){{\scriptsize $v$}}
 \end{picture}}
\def\BWc{
 \begin{picture}(22,10)(0,10)
 \multiput(4,10)(1.8,0){7}{\line(1,0){1.1}}
 \put(10,4){\line(0,1){12}}
 \put(10,10){\oval(2.5,2.5)[bl]}
 \put(2,9){{\scriptsize $j$}}
 \put(7,1.5){{\scriptsize $j-j^{\prime}$}}
 \put(16.5,9){{\scriptsize $j^{\prime}$}}
 \put(7,16.9){{\scriptsize $j^{\prime}-j$}}
 \put(7,7){{\scriptsize $v$}}
 \end{picture}}

\def\BWd{
 \begin{picture}(22,10)(0,10)
 \multiput(4,10)(1.8,0){7}{\line(1,0){1.1}}
 \put(10,4){\line(0,1){12}}
 \put(10,10){\oval(2.5,2.5)[br]}
 \put(2,9){{\scriptsize $j$}}
 \put(7,1.5){{\scriptsize $j^{\prime}-j$}}
 \put(16.5,9){{\scriptsize $j^{\prime}$}}
 \put(7,16.9){{\scriptsize $j-j^{\prime}$}}
 \put(11,7){{\scriptsize $v$}}
 \end{picture}}
\def\TM{
\begin{picture}(60,15)(0,15)
\put(0,10){\line(1,0){60}}
\put(0,20){\line(1,0){60}}
\put(10,5){\line(0,1){20}}
\put(20,5){\line(0,1){20}}
\put(40,5){\line(0,1){20}}
\put(50,5){\line(0,1){20}}
\put(10,10){\oval(2.5,2.5)[bl]}
\put(20,10){\oval(2.5,2.5)[bl]}
\put(40,10){\oval(2.5,2.5)[bl]}
\put(50,10){\oval(2.5,2.5)[bl]}
\put(10,20){\oval(2.5,2.5)[br]}
\put(20,20){\oval(2.5,2.5)[br]}
\put(40,20){\oval(2.5,2.5)[br]}
\put(50,20){\oval(2.5,2.5)[br]}
\put(1,7){{\scriptsize $u-iv$}}
\put(11,7){{\scriptsize $u-iv$}}
\put(31,7){{\scriptsize $u-iv$}}
\put(41,7){{\scriptsize $u-iv$}}
\put(11,17){{\scriptsize $u+iv$}}
\put(21,17){{\scriptsize $u+iv$}}
\put(41,17){{\scriptsize $u+iv$}}
\put(51,17){{\scriptsize $u+iv$}}
\put(8.5,1){{\footnotesize $V^{\prime}_1$}}
\put(18.5,1){{\footnotesize $V^{\prime}_2$}}
\put(38.5,1){{\footnotesize $V^{\prime}_{L-1}$}}
\put(48.5,1){{\footnotesize $V^{\prime}_{L}$}}
\put(-5,8.5){{\footnotesize $V_2$}}
\put(-5,17.5){{\footnotesize $V_{1}$}}
\end{picture}}

\def\QTM{
\begin{picture}(60,10)(0,10)
\put(0,10){\line(1,0){60}}
\put(10,5){\line(0,1){10}}
\put(20,5){\line(0,1){10}}
\put(40,5){\line(0,1){10}}
\put(50,5){\line(0,1){10}}
\put(10,10){\oval(2.5,2.5)[tr]}
\put(20,10){\oval(2.5,2.5)[br]}
\put(40,10){\oval(2.5,2.5)[tr]}
\put(50,10){\oval(2.5,2.5)[br]}
\put(11,12){{\scriptsize $u+iv$}}
\put(21,7){{\scriptsize $u-iv$}}
\put(41,12){{\scriptsize $u+iv$}}
\put(51,7){{\scriptsize $u-iv$}}
\put(8.5,1){{\footnotesize $V_1$}}
\put(18.5,1){{\footnotesize $V_2$}}
\put(38.5,1){{\footnotesize $V_{N-1}$}}
\put(48.5,1){{\footnotesize $V_{N}$}}
\put(-5,8.5){{\footnotesize $V^{\prime}_1$}}
\end{picture}}
\def\MQTM{
\begin{picture}(60,10)(0,10)
\put(0,10){\line(1,0){60}}
\put(10,5){\line(0,1){10}}
\put(20,5){\line(0,1){10}}
\put(40,5){\line(0,1){10}}
\put(50,5){\line(0,1){10}}
\put(10,10){\oval(2.5,2.5)[tr]}
\put(20,10){\oval(2.5,2.5)[br]}
\put(40,10){\oval(2.5,2.5)[tr]}
\put(50,10){\oval(2.5,2.5)[br]}
\put(11,12){{\scriptsize $u+iv$}}
\put(21,7){{\scriptsize $u-iv$}}
\put(41,12){{\scriptsize $u+iv$}}
\put(51,7){{\scriptsize $u-iv$}}
\put(8.5,1){{\footnotesize $V_1$}}
\put(18.5,1){{\footnotesize $V_2$}}
\put(38.5,1){{\footnotesize $V_{N-1}$}}
\put(48.5,1){{\footnotesize $V_{N}$}}
\put(-3,8.5){{\footnotesize $\ell$}}
\put(62,8.5){{\footnotesize $\ell^{\prime}$}}
\end{picture}}
\def\QTMF{
\begin{picture}(60,10)(0,10)
\multiput(0,10)(2,0){30}{\line(1,0){1}}
\put(10,5){\line(0,1){10}}
\put(20,5){\line(0,1){10}}
\put(40,5){\line(0,1){10}}
\put(50,5){\line(0,1){10}}
\put(10,10){\oval(2.5,2.5)[tr]}
\put(20,10){\oval(2.5,2.5)[br]}
\put(40,10){\oval(2.5,2.5)[tr]}
\put(50,10){\oval(2.5,2.5)[br]}
\put(11,12){{\scriptsize $u+iv$}}
\put(21,7){{\scriptsize $u-iv$}}
\put(41,12){{\scriptsize $u+iv$}}
\put(51,7){{\scriptsize $u-iv$}}
\put(8.5,1){{\footnotesize $V_1$}}
\put(18.5,1){{\footnotesize $V_2$}}
\put(38.5,1){{\footnotesize $V_{N-1}$}}
\put(48.5,1){{\footnotesize $V_{N}$}}
\end{picture}}

\def\Lattice{
\begin{picture}(90,45)(0,45)
\put(25,15){\line(0,1){60}}
\put(35,15){\line(0,1){60}}
\put(55,15){\line(0,1){60}}
\put(65,15){\line(0,1){60}}
\put(15,25){\line(1,0){60}}
\put(15,35){\line(1,0){60}}
\put(15,55){\line(1,0){60}}
\put(15,65){\line(1,0){60}}
\put(25,25){\oval(2.5,2.5)[bl]}
\put(35,25){\oval(2.5,2.5)[bl]}
\put(55,25){\oval(2.5,2.5)[bl]}
\put(65,25){\oval(2.5,2.5)[bl]}
\put(25,35){\oval(2.5,2.5)[br]}
\put(35,35){\oval(2.5,2.5)[br]}
\put(55,35){\oval(2.5,2.5)[br]}
\put(65,35){\oval(2.5,2.5)[br]}
\put(25,55){\oval(2.5,2.5)[bl]}
\put(35,55){\oval(2.5,2.5)[bl]}
\put(55,55){\oval(2.5,2.5)[bl]}
\put(65,55){\oval(2.5,2.5)[bl]}
\put(25,65){\oval(2.5,2.5)[br]}
\put(35,65){\oval(2.5,2.5)[br]}
\put(55,65){\oval(2.5,2.5)[br]}
\put(65,65){\oval(2.5,2.5)[br]}
\put(24,77){{\footnotesize $V_1^{\prime}$}}
\put(34,77){{\footnotesize $V_2^{\prime}$}}
\put(54,77){{\footnotesize $V_{L-1}^{\prime}$}}
\put(64,77){{\footnotesize $V_{L}^{\prime}$}}
\put(77,64){{\footnotesize $V_{1}$}}
\put(77,54){{\footnotesize $V_{2}$}}
\put(77,34){{\footnotesize $V_{N-1}$}}
\put(77,24){{\footnotesize $V_{N}$}}
\put(16,22){{\scriptsize $u-iv$}}
\put(26,22){{\scriptsize $u-iv$}}
\put(46,22){{\scriptsize $u-iv$}}
\put(56,22){{\scriptsize $u-iv$}}
\put(26,32){{\scriptsize $u+iv$}}
\put(36,32){{\scriptsize $u+iv$}}
\put(56,32){{\scriptsize $u+iv$}}
\put(66,32){{\scriptsize $u+iv$}}
\put(16,52){{\scriptsize $u-iv$}}
\put(26,52){{\scriptsize $u-iv$}}
\put(46,52){{\scriptsize $u-iv$}}
\put(56,52){{\scriptsize $u-iv$}}
\put(26,62){{\scriptsize $u+iv$}}
\put(36,62){{\scriptsize $u+iv$}}
\put(56,62){{\scriptsize $u+iv$}}
\put(66,62){{\scriptsize $u+iv$}}
\put(40,10){\vector(-1,0){25}}
\put(50,10){\vector(1,0){25}}
\put(10,50){\vector(0,1){25}}
\put(10,40){\vector(0,-1){25}}
\put(43,9){$L$}
\put(8,44){$N$}
\put(55,-5){\vector(0,1){12}}
\put(39,-2){$T_A(u,v)$}
\put(-8,40){\vector(1,0){12}}
\put(-8,44){$T_1(u,v)$}
\end{picture}}
\def\YBEa{\setlength{\unitlength}{0.6mm}
 \begin{picture}(50,25)(0,25)
 \put(10,5){\line(1,2){20}}
 \put(40,5){\line(-1,2){20}}
 \multiput(5,12.5)(2,0){20}{\line(1,0){1}}
 \bezier{50}(11,12.5)(11,11)(12.5,10)
 \bezier{50}(32.5,12.5)(33.5,8.5)(37.5,10)
 \bezier{50}(23.7,32.5)(25,31)(26.3,32.5)
 \put(8,8){{\scriptsize $v$}}
 \put(33,5){{\scriptsize $v^{\prime}$}}
 \put(27,32){{\scriptsize $v^{\prime}-v$}}
\end{picture}}
\def\YBEb{\setlength{\unitlength}{0.6mm}
 \begin{picture}(50,25)(0,25)
 \put(20,5){\line(1,2){20}}
 \put(30,5){\line(-1,2){20}}
 \multiput(5,37.5)(2,0){20}{\line(1,0){1}}
 \bezier{50}(33.5,37.5)(33.5,36)(35,35)
 \bezier{50}(10,37.5)(11,33.5)(15,35)
 \bezier{50}(23.7,12.5)(25,11)(26.3,12.5)
 \put(31,32){{\scriptsize $v$}}
 \put(8,30){{\scriptsize $v^{\prime}$}}
 \put(28,11){{\scriptsize $v^{\prime}-v$}}
\end{picture}}
\def\YBEc{\setlength{\unitlength}{0.6mm}
 \begin{picture}(50,25)(0,25)
 \put(10,5){\line(1,2){20}}
 \put(40,5){\line(-1,2){20}}
 \multiput(5,12.5)(2,0){20}{\line(1,0){1}}
 \bezier{50}(33.5,12.5)(33.5,13.5)(35,15)
 \bezier{50}(10,12.5)(11.5,16)(15,15)
 \bezier{50}(23.7,32.5)(25,31)(26.3,32.5)
 \put(8.5,14){{\scriptsize $v$}}
 \put(30,13){{\scriptsize $v^{\prime}$}}
 \put(27,32){{\scriptsize $v-v^{\prime}$}}
\end{picture}}
\def\YBEd{\setlength{\unitlength}{0.6mm}
 \begin{picture}(50,25)(0,25)
 \put(20,5){\line(1,2){20}}
 \put(30,5){\line(-1,2){20}}
 \multiput(5,37.5)(2,0){20}{\line(1,0){1}}
 \bezier{50}(11,37.5)(11,39)(12.5,40)
 \bezier{50}(32.5,37.5)(34,41)(37.5,40)
 \bezier{50}(23.7,12.5)(25,11)(26.3,12.5)
 \put(31,39){{\scriptsize $v$}}
 \put(7,38){{\scriptsize $v^{\prime}$}}
 \put(28,11){{\scriptsize $v-v^{\prime}$}}
\end{picture}}
\begin{document}
\title{Continued fraction TBA and functional relations\\
in XXZ model at root of unity} 
\author{Atsuo Kuniba\thanks{E-mail address:
atsuo@hep1.c.u-tokyo.ac.jp}, 
Kazumitsu Sakai\thanks{E-mail address:
sakai@as2.c.u-tokyo.ac.jp} and 
Junji Suzuki\thanks{E-mail address:
suz@hep1.c.u-tokyo.ac.jp} \\\\
\it Institute of Physics, University of Tokyo, \\
\it Komaba, Meguro-ku, Tokyo 153, Japan}
\maketitle
\begin{abstract}
Thermodynamics of the 
spin $\frac{1}{2}$ XXZ model is studied in the critical regime
using the quantum transfer matrix (QTM) approach.
We find functional relations indexed by the Takahashi-Suzuki numbers
among the fusion hierarchy of the
QTM's ($T$-system) and their certain 
combinations ($Y$-system).
By investigating analyticity of the latter,
we derive a closed set of non-linear integral equations which
characterize the free energy
and the correlation lengths for both 
$\langle \sigma_j^{+}\sigma_i^{-}\rangle$ and
$\langle \sigma_j^{z}\sigma_i^{z}\rangle$
at any finite temperatures.
Concerning the free energy, they
exactly coincide with Takahashi-Suzuki's
TBA equations based on the string hypothesis.
By solving the integral equations numerically 
the correlation lengths are determined,
which agrees with the earlier results
in the low temperature limit.

\end{abstract}
\section{Introduction}
In this paper we study the spin $\frac{1}{2}$
XXZ model at finite temperature based on the 
recently developed quantum transfer matrix (QTM) approach
\cite{MSuzPB}--\cite{JKSHub}.
We shall deal with the ``root of unity'' case in the 
gap-less regime.
Namely, the anisotropy parameter has the form 
$\Delta = \cos\!\frac{\pi}{p_0}$ with $p_0$ 
any rational number not less than 2.
(See (\ref{eqn:Hamiltonian}).)
We derive the non-linear integral equations 
that characterize the free energy
and the correlation lengths for both 
$\langle\sigma^+_j \sigma^-_i\rangle$ and 
$\langle\sigma^z_j \sigma^z_i\rangle$
at any finite temperatures.

Thermodynamics of the XXZ model is a classical 
and by no means fresh problem 
at least as far as the free energy is concerned.
It goes back to 1972 that Takahashi and Suzuki \cite{TS}
took the thermodynamic Bethe ansatz (TBA) approach \cite{YY}
to the free energy based on the elaborate string hypothesis.
They selected, as allowed lengths of 
strings, a special sequence of integers $n_j$
which we call the Takahashi-Suzuki (TS) numbers.
The resulting free energy yields
correct physical behaviours in many respects.
Actually this is one of the best known example among 
many successful applications
of the TBA and string hypotheses.
However there is also some 
controversy in string hypotheses themselves 
\cite{EKSstring,JuDo,AlczMart}, 
in view of which those successes are rather mysterious.

This is one of our motivations to revisit the XXZ model with 
the recent QTM method.
It integrates many ideas in the statistical mechanics
and solvable models \cite{MSuzPr}--\cite{JKSHub} and has a number of 
advantages over the traditional TBA approach.
It only relies on 
certain analyticity of the QTM, which can 
easily be confirmed much more convincingly by numerics.
Moreover it enables us to systematically calculate the correlation lengths
beyond the free energy for a wide range of temperatures.
See \cite{KZeit} for $\Delta > 1$ case.
Roughly, the QTM method goes as follows.
First one transforms the 1D quantum system into an integrable 
2D classical system based on the general equivalence theorem 
\cite{MSuzPr,MSuzPB}.
The QTM $T_1$ is a transfer matrix propagating in the 
cross channel of the latter.
Despite that the original 1D Hamiltonian $H$ is critical,
the QTM $T_1$ can be made to have a gap.
Therefore the formidable sum 
$\mbox{Tr} e^{-\beta H}$ ($\beta=1/k_B T$ : $T$ is 
temperature) can be expressed 
as its {\em single} eigenvalue which is largest in the magnitude.
Furthermore the correlation lengths are obtained from the ratio of 
the largest and the sub-leading eigenvalues of $T_1$.
To evaluate them actually, one must however recognize a price to pay; 
now the QTM $T_1(u)$ itself becomes dependent on the fictitious Trotter 
number $N$ through its size and also a coupling constant as
$T_1(u = \mbox{const}\frac{\beta}{N})$ \cite{SAW}--\cite{Mizu}.
This makes it difficult to determine the spectra of $T_1$ 
by a naive numerical extrapolation as $N \rightarrow \infty$.
A crucial idea to overcome this is to equip the QTM with
another variable $v$ and to 
exploit the Yang-Baxter integrability with respect to it;
$[T_1(u,v), T_1(u, v')] = 0$ \cite{Klu}.
Here $u$ and $v$ play the role of the (inverse) temperature and 
the spectral parameter, respectively.
Furthermore one introduces 
some auxiliary functions of $v$, which should realize 
somewhat miraculous features. 
Their appropriate combinations should have a nice analyticity 
that encodes the information on the Bethe ansatz roots of $T_1(u,v)$.
Once this is achieved, one can derive 
a non-linear integral equation 
which efficiently determines the 
sought eigenvalues of $T_1$.
The Trotter limit $N \rightarrow \infty$ 
can thereby be taken {\em analytically}.
The most essential step in this method is 
to invent such auxiliary functions and their appropriate combinations.
There are some interesting variety of choices for 
them in various models \cite{KZeit}--\cite{JKSHub}.

Back to the XXZ model our finding is that such 
auxiliary functions can be given by the QTM's 
$\{ T_{n_j-1} \vert \, n_j : \mbox{TS numbers}\}$, which is the subset 
of the known fusion hierarchy of commuting transfer matrices 
whose dimensions of the auxiliary spaces are 
precisely the TS numbers $n_j$.
We will show that 
$\{ T_{n_j-1} \}$ satisfy functional relations 
among themselves ($T$-system) and so do their certain ratios $\{Y_j\}$ 
an elaborate one ($Y$-system).
See (\ref{ysystem1})-(\ref{ysystem4}).
Especially there is a special 
identity (\ref{tsystem2}) among $\{T_{n_j-1}\}$ that holds only at
rational values of $p_0$ and makes the $Y$-system close 
finitely.
Besides the peculiarity at 
general roots of unity, 
such use of the fusion hierarchy as the auxiliary functions
originates in the studies of 
finite size corrections \cite{KP,KNS2}.

As for the free energy we thus obtain the  
integral equations identical with 
Takahashi-Suzuki's TBA equations but totally independently of 
their string hypothesis.
We shall further study the second
and the third largest eigenvalues of $T_1$ for $p_0$ integer.
They are related to the correlation length $\xi$ 
of $\langle\sigma^+_j \sigma^-_i \rangle$ and 
$\langle\sigma^z_j \sigma^z_i \rangle$, respectively.
In contrast with the largest eigenvalue, 
now the zeros of the $T_{n_j-1}$ come into the ``physical strip''
spoiling the nice analyticity.
Nevertheless we manage to identify their patterns and derive 
the ``excited state TBA equations''. 
Solving them numerically we determine the
curve $\xi = \xi(\beta)$.
Especially the low temperature asymptotics 
$\lim_{\beta \rightarrow \infty} \xi(\beta)/\beta$ agrees with the 
known result \cite{KZeit,KBIbook} with high accuracy
for the both correlations.

Our formulation here using the 2 variable 
QTM $T_1(u,v)$ and fusion hierarchies is based on 
\cite{Klu,JKSfusion}.
There are similar approaches 
in the context of integrable QFT's in a finite
volume \cite{Fen,BLZ,R}.

It has been known for some time that 
solutions of $Y$-systems can curiously be constructed 
from $T$-systems \cite{KP,KNS1}.
By now this connection has been generalized to arbitrary non-twisted 
affine algebra $X^{(1)}_r$ for the associated 
$Y$-system \cite{KN} and the $T$-system \cite{KNS1}.
(See also \cite{KS}.)
In this sense, our results here display a further connection 
of such sort for $U_q(\widehat{sl}(2))$ at $q$ general root of unity.

The layout of the paper is as follows. 
In section 2 we formulate the XXZ model at finite temperature
in terms of the QTM $T_1$.
In section 3 we give the fusion hierarchy $\{ T_{n-1} \}$ 
of QTM's and their eigenvalues.
A functional relation ($T$-system) valid for general $p_0$ 
is also given.
In section 4 we construct the $Y$-system out of the $T$-system.
The former closes finitely due to the special functional relation
(\ref{tsystem2}) valid only for rational $p_0$.
In section 5 we derive integral equations for the free energy, 
and in section 6 for 
the correlation lengths of $\langle\sigma^+_j \sigma^-_i \rangle$ and 
$\langle\sigma^z_j \sigma^z_i \rangle$.
Section 7 is a discussion.
Appendix A recalls the definition of the TS numbers and the related 
data.
Appendices B and C contain a check of the analyticity of the $Y$-functions for 
the free energy.
Appendix D is devoted to the free fermion case $p_0 = 2$, 
which needs a separate treatment.
\section{Quantum transfer matrix}
The Hamiltonian of spin $\frac{1}{2}$ 
one dimensional XXZ model on a periodic lattice with $L$ sites is
\bea
H&=&\sum_{j=1}^{L}H_{j\,j+1},\nonumber \\
H_{j\,j+1}&=&\frac{J}{4}\left(\sigma_{j}^{x}\sigma_{j+1}^{x}+
                            \sigma_{j}^{y}\sigma_{j+1}^{y}+
                         \Delta (\sigma_{j}^{z}\sigma_{j+1}^{z}-1)\right).
\label{eqn:Hamiltonian}
\eea
Here $\sigma_j^x,\sigma_j^y,\sigma_j^z$ are the local spin operators
(Pauli matrices) at the $j$-th lattice site and $J$ is a real 
coupling constant.
We shall consider the model with 
the anisotropy parameter $\Delta$ in 
the critical region $-1< \Delta < 1$.
Due to the invariance of the spectrum of (\ref{eqn:Hamiltonian})
under the transformation $(J,\Delta) \rightarrow (-J,-\Delta)$, 
we can further restrict the range to  $0\le\Delta<1$ and 
introduce the parametrization:
\bea
\Delta &=& \cos \theta \quad 0<\theta\le\frac{\pi}{2},\\
 p_0&=&\frac{\pi}{\theta} \quad p_0 \ge 2.
\eea
The model is associated with the 
quantum group $U_q(\widehat{sl}(2))$ 
at $q=e^{\pi i/p_0}$.

\begin{figure}[t]
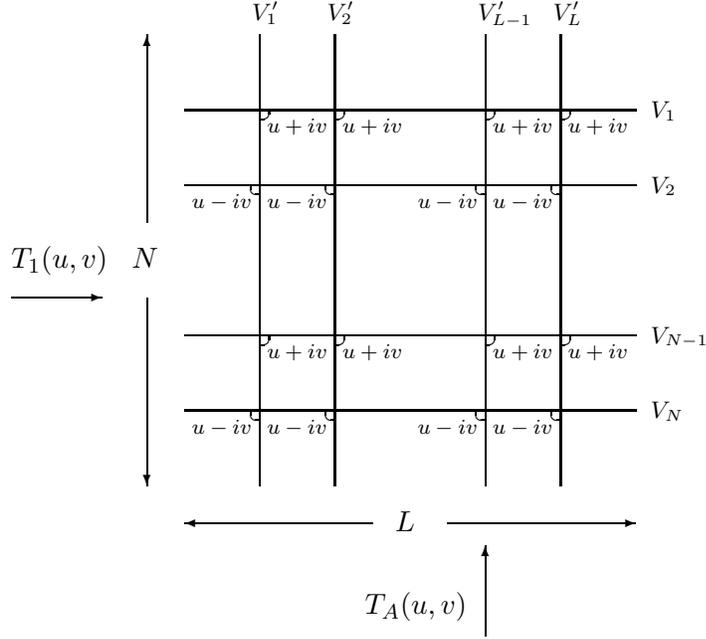

\begin{center} 
\Lattice \\[53mm]
\end{center}
\caption{Relation between $T_A(u,v)$ and $T_1(u,v).$}
\label{lattice}
\end{figure}

In order to consider its thermodynamics
we relate it to the six vertex model.
This is a 
two dimensional classical system whose Boltzmann weights are given by
\be
\vertex[\pm 1,\pm 1,\pm 1,\pm 1,v]=\frac{[v+2]}{[2]},\quad 
\vertex[\pm 1,\mp 1,\pm 1,\mp 1,v]=\frac{[v]}{[2]},\quad
\vertex[\pm 1,\mp 1,\mp 1,\pm 1,v]=1,\\[10mm]
\label{6vweights}
\ee
where
\be
[v]=\frac{\sin\frac{\theta}{2}v}{\sin\frac{\theta}{2}}.
\label{boxv}
\ee
Let $V$ be a two dimensional irreducible module over
$U_q(\widehat{sl}(2))$.
As is well known the quantum $R$-matrix 
$R \in \mbox{End}(V \otimes V)$ with the above matrix elements 
and the spectral parameter $v$ satisfies 
the Yang-Baxter equation (YBE) (cf.(\ref{eqn:YBE})).
To relate the six vertex model with the XXZ model, 
consider a two dimensional square 
lattice with $N$ rows and $L$ columns.
We shall assume that 
$N$ is even throughout.
Define the (auxiliary) transfer matrix $T_A(u,v)$ as
\be
T_A(u,v)=\mbox{Tr}_{V_{1} \otimes V_2}
\Biggl(\quad \quad \TM \quad \Biggr).\\[15mm]
\label{eqn:TM}
\ee
See also Fig.\ref{lattice}.
Here and in what follows $V_1 = \cdots = V_N = V'_1 = \cdots = V'_L = V$.
Operators diagrammatically shown as in (\ref{eqn:TM}) are always 
assumed to act on the states in the bottom line to transfer them into those
in the upper line.
Using the identity 
($P_{j\,j+1}$ is the permutation operator acting on
$V_j \otimes V_{j+1}$),
\bean
R_{j\,j+1}(0)\frac{d}{dv}R_{j\,j+1}(v)\bigg|_{v=0}&=&
P_{j\,j+1}   \frac{d}{dv}R_{j\,j+1}(v)\bigg|_{v=0} \\
&=& \frac{\theta}{J \sin \theta}
\left(H_{j\,j+1} +\frac{J}{2}\Delta\right),
\eean
we expand $T_A(u,0)$ as 
\be
T_A(u,0)=1+\frac{2\theta u}{J\sin \theta}\left(H+\frac{JL \Delta}2
\right)+O(u^2).
\label{eqn:baxter}
\ee
This formula represents 
as an equivalence of the XXZ model and the six vertex model.
In fact we can go further to the finite temperature case.
{}From (\ref{eqn:baxter}) it follows that
\be
\exp\left(-\beta\left(H+\frac{JL}{2}\Delta\right)\right)=
\lim_{N \to \infty}T_A(u_N,0)^{\frac{N}{2}}, 
\quad \quad u_N=-\frac{\beta J\sin\theta}{\theta N}. \label{uenu}
\ee 
Thus the free energy per site $f$ of the XXZ model is given by
\be
f=-\lim_{L \to \infty}\lim_{N \to \infty}\frac{1}{L\beta}
   \ln \left(\mbox{Tr}_{V^{\otimes L}}T_{A}(u_N,0)^{\frac{N}{2}} \right)-
   \frac{J}{2} \Delta.
\ee
However, eigenvalues of the transfer matrix $T_A(u_N,0)$ 
are infinitely degenerate in the limit $u_N \stackrel{N \to \infty}
{\longrightarrow} 0$, therefore the 
the trace renders a serious problem.
To avoid this,
we introduce the following trick; we rotate the lattice by $90^{\circ}$ 
(see Fig.\ref{lattice}) and rewrite the free energy as
\be
f=-\lim_{L \to \infty}\lim_{N \to \infty}\frac{1}{L\beta}
\ln \left( \mbox{Tr}_{V^{\otimes N}}T_1(u_N,0)^{L} \right)-
\frac{J}{2} \Delta,
\ee
where
\be
T_1(u,v)=\mbox{Tr}_{V^{\prime}_1}
\Biggl(\quad \quad \QTM \quad \Biggr).\\[13mm]
\label{eqn:QTM}
\ee
We call $T_1(u,v)$ the quantum transfer matrix (QTM).
Note that due to the YBE (\ref{eqn:YBE}), the QTM is
commutative as long as the 
$u$ variable is taken same:
\[
 [T_1(u,v),T_1(u,v^{\prime})]=0.
\]
{}From now on, we write 
the $k$-th largest eigenvalue of the matrix
$T_1(u,v)$ as $T^{(k)}_1(u,v)$.
Since the two limits are exchangeable as proved in \cite{MSuzPB,InSuz},
we take the limit $L \to \infty$ first.
Noting that there is a finite gap between
$\lim_{N \to \infty} T^{(1)}_1(u_N,0)$ and
$\lim_{N \to \infty} T^{(2)}_1(u_N,0)$, we have 
\be
f=-\frac{1}{\beta} 
\lim_{N \to \infty} \ln T_1^{(1)}(u_N,0)-\frac{J}{2}\Delta.
\label{eqn:ebaxter}
\ee
Namely, the problem of describing the thermodynamics of one dimensional 
quantum systems reduces to finding the largest eigenvalue of the QTM of two
dimensional finite systems (to be exact, finite in 
the vertical direction only).
Equation (\ref{eqn:ebaxter}) is a 
finite temperature extension of the equation (\ref{eqn:baxter}).

In this approach, the thermodynamical completeness
$-\lim_{\beta \rightarrow 0}\beta f = \ln 2$ follows 
easily from  $T^{(1)}_1(0,0) = 2$, which is
obvious from $R_{12}(0)=P_{12}$.
(As for the combinatorial completeness see 
\cite{KL} including the higher spin cases.)

Most significantly this method makes it possible 
to calculate some correlation length $\xi_k$
($k\ge 2$) at finite temperature. 
To see this 
let $Q_i = 1 \otimes \cdots \otimes Q \otimes  \cdots 
\otimes 1 \in \hbox{End}(V^{\otimes L})$ be a local operator acting on the 
$i$-th site $V'_i$ via $Q = \sum_{\ell,\ell' = \pm 1}
Q_{\ell', \ell}E_{\ell', \ell}$.
Here $E_{\ell', \ell}$ denotes the $2$ by $2$ elementary matrix 
and $Q_{\ell', \ell}$
is the matrix element.
($\sigma^z = E_{1,1}-E_{-1,-1}, 
\sigma^+ = E_{1,-1}, \sigma^- = E_{-1,1}$.)
Given $Q$ we introduce the operator 
$S(Q \vert u,v) \in \hbox{End}(V^{\otimes N})$ by
$$
S(Q \vert u,v)=\sum_{\ell,\ell^{\prime}=\pm 1}Q_{\ell^{\prime},\ell}
              \Biggl(\quad \MQTM \quad \Biggr).\\[13mm]
$$
Then for the local operators $Q_i, P_j \in \hbox{End}(V^{\otimes L})$ their
finite temperature correlation function is expressed as 
\bean
 \langle P_j Q_i \rangle &=&
    \lim_{L \rightarrow \infty}
  \frac{\mb{Tr}_{V^{\otimes L}}\,P_j Q_i \exp(-\beta H)}
     {\mb{Tr}_{V^{\otimes L}}\,\exp(-\beta H)}  \\
&=& \lim_{N\to\infty}\lim_{L \to \infty}
    \frac{\mb{Tr}_{V^{\otimes N}}\left(S(P \vert u_N,0)
    T_1(u_N,0)^{j-i-1}S(Q \vert u_N,0) T_1(u_N,0)^{L+i-j-1}\right)}
    {\mb{Tr}_{V^{\otimes N}} T_1(u_N,0)^L},
\eean
where we have exchanged the two limits.
Suppose that  $Q^{(k)}$ and $Q^{(k)\dagger}$ 
are the operators such that 
the matrix elements 
of $S(Q^{(k)} \vert u_N,0)$ and $S(Q^{(k)\dagger} \vert u_N,0)$ 
between the eigenspaces for
$T^{(1)}_1(u_N,0)$ and $T^{(k^{\prime})}_1(u_N,0)$
are zero for $1\le k^{\prime}<k$ and nonzero for $k^{\prime}=k$.
Then setting $(P,Q)=(Q^{(k)\dagger},Q^{(k)})$ in the
above, we have 
$$
\langle Q^{(k)\dag}_jQ^{(k)}_i \rangle \sim 
\lim_{N\to\infty}\left(\frac{T^{(k)}_1(u_N,0)}
{T^{(1)}_1(u_N,0)}\right)^{j-i}\qquad  k \ge 2,\, j \gg i.
$$
Fitting this with $\exp(-(j-i)/\xi_k)$ in the limit $j\gg i$
we have
\be
\frac{1}{\xi_k}=-\lim_{N \to \infty}\ln \left|\frac{T^{(k)}_1(u_N,0)}
{T^{(1)}_1(u_N,0)}
                \right|.
\label{cl}
\ee
As seen in section 6, $\xi_2$ and $\xi_3$ are 
the correlation lengths of 
$\langle\sigma_j^{+}\sigma_i^{-}\rangle$ ($Q^{(2)}=\sigma^{-}$)
and
$\langle\sigma_j^{z}\sigma_i^{z}\rangle$ ($Q^{(3)}=\sigma^{z}$),
respectively. 
\section{$T$-system}

To study $T_1(u,v)$, we embed it into a more 
general family of transfer matrices and 
explore the functional relations that govern 
the total system.
For this we consider the fusion hierarchy
$\{ T_{n-1}(u,v) \}$ defined by 
\be
T_{n-1}(u,v) = \mbox{Tr}\Biggl( \QTMF \Biggr), \\[10mm]
\label{ftdefinition}
\ee
where the trace is over the $n$ dimensional 
irreducible auxiliary space depicted by the dotted line.
To be explicit, we give the constituent fusion 
Boltzmann weights.
\begin{eqnarray}\label{fbw}
&&\BWa =  \BWb = \frac{[v+1 \pm (n+1-2j)]}{[2]},\nonumber\\[10mm]
&&\BWc = \epsilon_{j j'} 
\frac{\sqrt{[2\mbox{min}(j,j')][2n-2\mbox{min}(j,j')]}}{[2]},
\nonumber\\[10mm]
&&\BWd = \epsilon'_{j j'} 
\frac{\sqrt{[2\mbox{min}(j,j')][2n-2\mbox{min}(j,j')]}}{[2]}.
\nonumber\\[10mm]
\end{eqnarray}
Here $j, j' \in \{ 1, \ldots, n \}$ and 
$\vert j - j' \vert = 1$.
$\epsilon_{j j'}, \epsilon'_{j j'}$ are arbitrary
parameters such that 
$\epsilon_{j j'}\epsilon_{j' j} 
= \epsilon'_{j j'}\epsilon'_{j' j} = 1$.
If they are 1, the six vertex case $n=2$ of these 
weights reduce to (\ref{6vweights}) 
under the identification of $j=1$ and $2$ states 
with $+1$ and $-1$, respectively.
The Boltzmann weights satisfy the YBE:
\be
\YBEa=\YBEb, \YBEc=\YBEd \\[15mm]
\label{eqn:YBE}
\ee
{}From the picture (\ref{ftdefinition}) one sees that 
the members of the fusion 
hierarchy are all commutative for the same $u$:
$$
[ T_{n-1}(u,v), T_{n'-1}(u,v')] = 0,
$$
due to the $R$-matrix intertwining the $n$ and $n'$ 
dimensional representations.
Thus they can be simultaneously diagonalized and the eigenvalues 
(also written as $T_{n-1}(u,v)$) 
are readily obtained in the dressed vacuum form:
\begin{eqnarray}
&T_{n-1}(u,v) &= \sum_{j=1}^n
\phi\bigl(v-i(u+n+2-2j)\bigr)
\phi\bigl(v+i(u-n+2j)\bigr) \nonumber \\
&&\times\frac{Q(v+in)Q(v-in)}
{Q\bigl(v+i(2j-n)\bigr)Q\bigl(v+i(2j-n-2)\bigr)},
\label{dvf}\\
&\phi(v) &= 
\biggl(\frac{\sh \frac{\theta}{2}v }{\sin \theta}\biggr)^{\frac{N}{2}},
\label{phidefinition}\\
&Q(v) &= \prod_{j=1}^m \sh \frac{\theta}{2}(v-\omega_j).
\label{qdefinition}
\end{eqnarray}
Here $m \in \{0, 1, \ldots, N/2 \}$ is the quantum number
counting the $(-1)$-states on odd sites and 
$(+1)$-states on even sites.
The dressed vacuum form is built upon the pseudo vacuum state 
$\bigl((+1) \otimes (-1)\bigr)^{\otimes \frac{N}{2}}$, 
which corresponds to $m=0$.
$\{ \omega_j \}$ is a solution of the Bethe ansatz equation 
(BAE):
\begin{equation}
-\Biggl( 
\frac{\sh\frac{\theta}{2}\bigl(\omega_j + i(u+2)\bigr)
\sh\frac{\theta}{2}(\omega_j - iu)}
{\sh\frac{\theta}{2}\bigl(\omega_j - i(u+2)\bigr)
\sh\frac{\theta}{2}(\omega_j + iu)}\Biggr)^{\frac{N}{2}}
= \frac{Q(\omega_j + 2i)}{Q(\omega_j - 2i)}.
\label{bae}
\end{equation}
The largest eigenvalue of $T_1(u,v)$ lies in 
the sector $m = N/2$.
Note that 
$T_{-1}(u,v) = 0$ and 
$T_0(u,v) = \phi(v-i(u+1))\phi(v+i(u+1))$.
An important property is the periodicity
\begin{equation}
T_{n-1}(u,v) = T_{n-1}(u,v+2p_0 i).
\label{tperiod}
\end{equation}

Let us present the functional relations among the
fusion hierarchy.
For any $v \in {\bf C}$ and integers $n \ge y \ge 1$,
the following is valid, which we 
call the $T$-system.
\begin{equation}
T_{n-1}(v + i y)T_{n-1}(v - i y) = 
T_{n+y-1}(v)T_{n-y-1}(v) +
T_{y-1}(v + i n)T_{y-1}(v - i n).
\label{tsystem1}
\end{equation}
Hereafter we shall often omit the common $u$ 
variable to simplify the notation.
The proof of this equation is direct by using the expression
(\ref{dvf}).
Representation theoretically, it 
is a simple consequence of the 
general exact sequence in \cite{CP}
as explained in \cite{KNS1} for $y=1$.
In general the $T$-system (\ref{tsystem1}) extends over 
infinitely many transfer matrices.
However, as we shall see in 
the next section, 
for rational $p_0$ there is a special functional relation (\ref{tsystem2})
that makes the associated $Y$-system closes finitely.

\section{$Y$-system at Root of Unity}

{}From now on, we shall concentrate on the case
when $p_0 > 2$ is a rational number and treat the free fermion case
$p_0$=2 separately in Appendix D.
Consider the continued fraction expansion of 
$p_0$ 
\begin{equation}
p_0 = \nu_1 + \frac{1}
      {\displaystyle \nu_2 + \frac{1}
        {\displaystyle \frac{\ddots}
         {\displaystyle \nu_{\alpha-1} + \frac{1}
          {\displaystyle \nu_\alpha}}}},
\label{p0exp}
\end{equation}
which 
specifies $\alpha \ge 1$ and 
$\nu_1, \ldots,  \nu_\alpha \in {\bf Z}_{\ge 1}$.
{}From the assumption $p_0 > 2$, we have $\nu_1 \ge 2$.
In fact $\nu_1 = 2$ is allowed only if $\alpha \ge 2$, and
$p_0 = \nu_1 \ge 3$ is assumed if $\alpha =1$.

In Appendix A
we recall the sequences of numbers
$\{m_j\}_{j=0}^{\alpha+1},
\{p_j\}_{j=0}^{\alpha+1}, 
\{y_j\}_{j=-1}^{\alpha}, 
\{z_j\}_{j=-1}^{\alpha}$ and $\{n_j\}_{j\ge 1}$ 
introduced in \cite{TS}. 
The last one is the TS numbers.
We shall also introduce its slight rearrangement 
$\{\widetilde{ n}_j\}_{j\ge 1}$ and a 
similar sequence $\{w_j\}_{j\ge 1}$ related to the 
``parity'' of the TS strings.
They are all specified uniquely from $p_0$. 
With those definitions we now describe a 
functional relation of the $T_{n-1}(v)$,
which is valid only at the 
root of unity and is
relevant to our subsequent argument.
\begin{equation}
T_{y_\alpha + y_{\alpha-1}-1}(v) = 
T_{y_\alpha - y_{\alpha-1}-1}(v) + 
2(-1)^{m z_\alpha}T_{y_{\alpha-1}-1}(v+iy_\alpha).
\label{tsystem2}
\end{equation}
Here $m$ is the number of the BAE roots in (\ref{qdefinition}).
The proof is straightforward by using 
(\ref{dvf}), (\ref{yzproperty2}) and 
$Q(v+2iy_\alpha) = (-1)^{m z_\alpha}Q(v)$.
When $\alpha=1$ hence  $p_0 = \nu_1$, (\ref{tsystem2}) 
reduces to a simple relation 
$T_{\nu_1}(v) = 
T_{\nu_1-2}(v) + 
2(-1)^m T_0(v+i\nu_1)$.

For $1 \le j \le j_{max} := m_\alpha -1$, 
let $0 \le r \le \alpha - 1$ be the unique integer
satisfying $m_r \le j < m_{r+1}$. 
Set 
\begin{eqnarray}
&Y_j(v) &= \frac{T_{{\tilde n}_{j+1}+y_r -1}(v + iw_jp_0)
               T_{ {\tilde n}_{j+1}-y_r -1}(v + iw_jp_0)}
              {T_{y_r -1}(v + i{\widetilde n}_{j+1} + iw_jp_0)
               T_{y_r -1}(v - i{\widetilde n}_{j+1} + iw_jp_0)},
\label{ydef1}\\
&1 + Y_j(v) &= \frac{T_{{\tilde n}_{j+1} -1}(v + iy_r + iw_jp_0)
               T_{{\tilde n}_{j+1} -1}(v -iy_r + iw_jp_0)}
              {T_{y_r -1}(v + i{\widetilde n}_{j+1} + iw_jp_0)
               T_{y_r -1}(v - i{\widetilde n}_{j+1} + iw_jp_0)},
\label{ydef2}\\
&K(v) &= \frac{(-1)^{m z_\alpha}T_{{\tilde n}_{j_{max}} -1}
                  (v + iw_{j_{max}}p_0)}
                 {T_{y_{\alpha-1} -1}
                  (v + iy_\alpha + iw_{j_{max}}p_0)},
\label{ydef3}
\end{eqnarray}
where in (\ref{ydef3}) 
${\widetilde n}_{j_{max}} = y_\alpha - y_{\alpha-1}$ and 
$w_{j_{max}} = z_\alpha - z_{\alpha-1} - 1$ in accordance with 
(\ref{ndefinition}) and (\ref{wdefinition}).
We also set $Y_0(v) = 0$ and $Y_{-1}(v) = \infty$.
Thanks to the $T$-system (\ref{tsystem1}), 
(\ref{ydef1}) and (\ref{ydef2}) are equivalent.
We find that $\{Y_j(v) \}_{j=1}^{m_\alpha - 1}$ and $K(v)$
close among the following finite set of functional relations,
which we call the $Y$-system.
\begin{th}
\begin{eqnarray}
&\hbox{For }& m_{r-1} \le j \le m_r - 2\,\, ( 1 \le r \le \alpha), \nonumber \\
&&Y_j(v+ip_r)Y_j(v-ip_r) = \bigl( 1+Y_{j-1}(v)\bigr)^{1-2\delta_{j m_{r-1}}}
\bigl(1+Y_{j+1}(v)\bigr),\label{ysystem1}\\
&\hbox{for }& j = m_r -1\,\, (1 \le r \le \alpha -1 ),  \nonumber \\
&&Y_j(v + ip_r + ip_{r+1})Y_j(v + ip_r - ip_{r+1})
Y_j(v - ip_r + ip_{r+1})Y_j(v - ip_r - ip_{r+1}) \nonumber \\
&&= \biggl\{ \bigl(1+Y_{j-1}(v+ip_{r+1})\bigr)
\bigl(1+Y_{j-1}(v-ip_{r+1})\bigr) \biggr\}^{1-2\delta_{1 \nu_r}}
\bigl(1+Y_{j+1}(v+ip_r)\bigr)\nonumber \\
&&\times \bigl(1+Y_{j+1}(v-ip_r)\bigr)
\bigl(1+Y_j(v + ip_r - ip_{r+1})\bigr)
\bigl(1+Y_j(v - ip_r + ip_{r+1})\bigr),
\label{ysystem2}\\
&&1 + Y_{m_\alpha - 1}(v) = \left(1 + K(v) \right)^2,
\label{ysystem3}\\
&&K(v + ip_\alpha)K(v - ip_\alpha) = 
1 + Y_{m_\alpha - 2}(v).
\label{ysystem4}
\end{eqnarray}
\end{th}
This can be proved by combining the $T$-systems 
(\ref{tsystem1}), (\ref{tsystem2}) with 
the definitions of 
$\{m_j\}, \{p_j\}, \{y_j\}, 
\{\widetilde{ n}_j\}$ and $\{w_j\}$ 
in Appendix A.
When $\nu_r = 1$, (\ref{ysystem1}) is void and 
(\ref{ysystem2}) holds for $j = m_r -1 = m_{r-1}$.
See (\ref{msec}).

In the case $p_0\in {\bf Z}_{\ge3}$, the $Y$-system
has a simple form  ($1\le j \le p_0-2$)
\bea
&&Y_{j}(v+i)Y_{j}(v-i)=(1+Y_{j-1}(v))(1+Y_{j+1}(v)),\label{ysystemc1} \\
&&1+Y_{p_0-1}(v)=(1+K(v))^2,\label{ysystemc2} \\
&&K(v+i)K(v-i)=1+Y_{p_0-2}(v),\label{ysystemc3}
\eea
where 
\bea
&&Y_{j}(v)=\frac{T_{j+1}(v)T_{j-1}(v)}{T_0(v+i(j+1))T_0(v-i(j+1))},
\label{yfnc1} \\
&&1+Y_{j}(v)=\frac{T_{j}(v+i)T_{j}(v-i)}
             {T_0(v+i(j+1))T_0(v-i(j+1))},\label{yfnc2} \\
&&K(v)=(-)^m \frac{T_{p_0-2}(v)}{T_0(v+ip_0)}
\label{yfnc3}
\eea
for $1\le j \le p_0-2$.
Due to the property of the TS-number, $Y_1(v)$ and $1+Y_1(v)$
are always given by setting $j=1$ in 
(\ref{yfnc1}) and (\ref{yfnc2}) for arbitrary $p_0>2$. 

A similar 
$Y$-system has been considered in \cite{Ta} in
a perturbed CFT context.
\section{Integral equation for free energy}
Let $\vert k \rangle$ be the eigenvector corresponding to the 
$k$-th largest eigenvalue $T^{(k)}_1(u,v)$ of $T_1(u,v)$.
We define the $k$-th (not necessarily $k$-th largest) 
eigenvalue $T^{(k)}_{n-1}(u,v)$ of the auxiliary QTM 
$T_{n-1}(u,v)$ by 
$T_{n-1}(u,v) \vert k \rangle  =  T^{(k)}_{n-1}(u,v) \vert k \rangle$.
Let $\{Y^{(k)}_j\}$ (and $K^{(k)}$) be the $Y$-functions
constructed from $\{T^{(k)}_{n-1}\}$ as in 
(\ref{ydef1}) -- (\ref{ydef3}).
In this section, we study the analyticity of 
$\{Y^{(1)}_j(u,v)\}$ and $K^{(1)}(u,v)$ in the complex $v$-plane.
Then we derive the integral equations which 
characterize the free energy.

An advantage in the present approach lies in the
fact  that the analytic assumption given below can be
explicitly checked numerically keeping the Trotter number $N$ finite.
We have performed  numerical studies with various values of $p_0$, 
$\beta$ and $N$ in determining the location 
zeros of fusion QTM's.
For example, the zeros for $T^{(1)}_{n-1}(u,v)$ for 
$p_0= \frac{24}{5}, u=\pm 0.01 $, $n$=2,3,4,5 $N=16$ and $N=32$ are plotted in
Fig \ref{zeros1}.
Guided by them we have the following for $u$ small 
(typically $\vert u \vert \sim 0.01$).
\begin{cj}
All the zeros of $T^{(1)}_{n-1}(u,v)$ are located on the
 almost straight lines 
$\Im v =  \pm n  \hbox{ mod } 2p_0$.
\end{cj}
This coincides with the observation in the XXX model
if one forgets the periodicity  $2p_0$
in the imaginary direction.
It is also consistent with the functional relation 
(\ref{tsystem1}).
The deviation from the straight line 
is very small ($\sim 10^{-2}$ at most) as  seen in the figures.
It becomes smaller as $u \rightarrow  0$.
Once Conjecture 1 is assumed, we can identify the strips
in the complex $v$-plane in which
$Y^{(1)}_j(v)$ or $ 1+Y_j^{(1)}(v)$ are {\em a}nalytic, 
{\em n}on{\em z}ero and
have {\em c}onstant asymptotics at $v = \pm \infty$. 
We call this property ANZC.
In Appendix B we verify that $Y^{(1)}_{j}(v)$, for example, is 
ANZC in the strip $\vert \Im v \vert \le x$ whenever the 
combination $Y^{(1)}_{j}(v+ix)Y^{(1)}_{j}(v-ix)$ takes place
in the $Y$-system (\ref{ysystem1})--(\ref{ysystem4}).
Apart from the exceptional 
Case 1, 2 and 3 listed below,
this makes it possible to transform 
most of the $Y$-system 
into integral equations defined on the real axis quite easily.
This is a consequence of a simple lemma.
To present it we let 
$\st[x]$ denote the strip $\Im v \in[-x,x]$ 
in the complex $v$-plane ($x \in {\bf R}_{>0}$).
Then we have
\begin{lem}
Suppose the functions $g_i(v)$ satisfy
\begin{equation}
g_0(v-iv_0)g_0(v+ iv_0) = 
\prod_{j\ge1}g_j(v-iv_j) g_j(v+iv_j),
\end{equation}
where $v_j \ge 0$ are real numbers and $v_0 >v_j (j \ge 1)$.
Assume further that $g_j(v)$ is ANZC 
in the strip $\st[ w_j]$ for some $w_j \ge v_j$ for $j \ge 0$. 
Then the above functional relation can be transformed into
the integral equation
\begin{eqnarray*}
 \ln g_0(v) &=& \sum_{j\ge 1} \int_{-\infty}^{\infty} 			
    R_j(v-v') \ln g_j(v')  dv' +\hbox{constant},  \\
    R_j(v) &=& \frac{1}{2 \pi} \int_{-\infty}^{\infty}
             e^{ikv} \frac{\ch v_j k}{\ch v_0 k} dk,
\end{eqnarray*}
where the constant is determined by the asymptotic values of the both sides.
\end{lem}

The proof uses Cauchy's theorem and the fact that the ANZC function 
$g_j(v)$ admits the Fourier transformation of its 
logarithmic derivative.

There are few exceptions to which the above lemma can not be applied
directly:
\begin{itemize}
\item  Case 1. $j=1$ in (\ref{ysystem1}),
\item  Case 2. $r=1$ when $\nu_1= 3, \nu_2=1$ and $\nu_3 \gg 1$ 
in (\ref{ysystem2}),
\item  Case 3. $r=1$ when $\nu_1=2$ in (\ref{ysystem2}).
\end{itemize}
Nevertheless, they can still be converted into integral equations after 
a suitable recipe.
Let us explain this for the most important Case 1 below.

Case 1 in (\ref{ysystem1}) is explicitly given by
\begin{equation}
Y^{(1)}_1(v+i) Y^{(1)}_1(v-i) = 1+Y^{(1)}_2(v).
\label{Ysys1}
\end{equation}
$Y^{(1)}_1(v)$ possesses zeros (resp. poles) of order 
$N/2$ at $\pm(1+u)i$ (resp. $\pm(1-u)i$)  mod $2p_0i$
in the strip  $\st[1]$ for $J > 0$ (resp. $J < 0$).
(Note that $u = u_N$ is a small quantity given in (\ref{uenu}.) 
Thus the lhs of (\ref{Ysys1}) does not meet the
condition for Lemma 1.
A simple trick, however, makes it applicable.
Define a modified function
\begin{equation}
\widetilde{Y}_1^{(1)}(u,v) = \frac{Y^{(1)}_1(u,v)}{
  ( \tnh\frac{\pi}{4}(v-i(1 \pm u))\tnh\frac{\pi}{4}(v+i(1 \pm u))
)^{\pm N/2}},
\end{equation}
where the $+$ and $-$ signs in front of $u$ and 
$N$ are chosen according as $J > 0$ and $J < 0$, respectively.
Then $\widetilde{Y}^{(1)}_1(v)$ has the ANZC property in $\st[1]$.
Due to the trivial identity $\tnh\frac{\pi}{4}(x+i) \tnh\frac{\pi}{4}(x-i)=1$,
$Y^{(1)}_1(v\pm i)$ in the lhs of (\ref{Ysys1}) can be replaced by 
$\widetilde{Y}^{(1)}_1(v\pm i)$.
Now the lemma applies.  The asymptotic values of both sides can be
immediately evaluated from the explicit results on the $T$-functions.
Then we have,
\begin{eqnarray}
 \ln Y^{(1)}_1(u,v) &=& 
    \pm \frac{N}{2} 
\ln  ( \tnh\frac{\pi}{4}(v-i(1 \pm u))\tnh\frac{\pi}{4}(v+i(1 \pm u))) 
	     \nonumber \\
    &+&
    \int_{-\infty}^{\infty}  \frac{1}{ 4 \ch\frac{\pi (v-v')}{2}   }
    \ln(1+Y_2^{(1)})(v')  dv'.
\label{Yint1}
\end{eqnarray}
Cases 2 and 3 are discussed in Appendix C.
In this way all the $Y$-system can be transformed into
coupled integral equations.
For finite $N$ one can evaluate $Y^{(1)}_j$'s given by
(\ref{ydef1})--(\ref{ydef3}) and 
(\ref{dvf})--(\ref{qdefinition}) directly from the BAE roots.
Or one can solve the integral equations numerically.
We have checked that the two independent calculations
lead to the same result up to $N=40$.

Let us proceed to the Trotter limit $N \rightarrow \infty$.
{}From now on we write the $Y$-functions in the limit as
\bea
\eta^{(k)}_j(v)&=&\lim_{N\to\infty}Y_j^{(k)}(u_N,v),
\label{etadef}\\
\kappa^{(k)}(v)&=&\lim_{N\to\infty}K^{(k)}(u_N,v).
\label{kappadef}
\eea
Apart from the $Y$-functions the $N$-dependence enters 
(\ref{Yint1}) only through the ``driving term''. Its 
large $N$ limit can be taken analytically as 
\[
\pm  \lim_{N\to\infty}\frac{N}{2}\ln \biggl( \tnh\frac{\pi}{4}(v-i(1 \pm u_N))
                                 \tnh\frac{\pi}{4}(v+i(1 \pm u_N)) \biggr)
    =-\frac{\beta\pi J \sin\theta}{2\theta\ch(\pi v/2)}.
\]
We thus arrive at the integral equations for 
$\eta^{(1)}_j$ and $\kappa^{(1)}$ which are independent of 
the fictitious Trotter number $N$.
\begin{eqnarray}
\ln \eta^{(1)}_j(v)&=&
       -\frac{\beta \pi J \sin\theta}{2 \theta \ch(\pi v/2)} \delta_{j,1} +
          (1-2 \delta_{m_{r-1},j}) s_r * \ln(1+\eta^{(1)}_{j-1})(v) 
                                                       \nonumber \\
          &+& s_r*\ln(1+\eta^{(1)}_{j+1})(v)
          \qquad \hbox{for } m_{r-1}\le j \le m_r-2, j\ge 1, 1\le r \le 
\alpha, 
		  \label{nlie1} \\
\ln \eta^{(1)}_j(v)&=&
 -\frac{\beta \pi J \sin\theta}{2 \theta \ch(\pi v/2)} \delta_{j,1} +
  (1-2 \delta_{m_{r-1},j}) s_r * \ln(1+\eta^{(1)}_{j-1})(v) +
             d_r*\ln(1+\eta^{(1)}_j)(v)                  \nonumber  \\
          &+& s_{r+1}*\ln(1+\eta^{(1)}_{j+1})(v) 
          \qquad \hbox{for } j=m_r-1, 1\le r <  \alpha, \label{nlie2}\\
\ln \kappa^{(k)}(v) &=& s_{\alpha}*\ln (1+\eta^{(1)}_{m_{\alpha}-2})(v),
\label{nlie3}
\end{eqnarray}
where $A*B(v)$ denotes the convolution 
$\int^\infty_{-\infty} A(v-v') B(v') dv'$, and
\be
 s_r(v)=\frac{1}{4p_r \ch \frac{\pi v}{2p_r}},  \qquad
       d_r(v)=\int^\infty_{-\infty} e^{ikv} 
              \frac{\ch(p_{r}-p_{r+1})k }{4\pi \ch (p_r k) \ch (p_{r+1} k)}
              dk.  
\label{nlie4}
\ee
The set of the equations closes by one further algebraic equation:
$\eta^{(1)}_{m_{\alpha}-1}(v)=\kappa^{(1)}(v)^{2}+2\kappa^{(1)}(v)$.
Under the identification $\eta^{(1)}_j(v) = \eta_j(v)$ and
$\kappa^{(1)}(v)=\kappa(v)$, the eqs. (\ref{nlie1})--(\ref{nlie3})
are nothing but the TBA equation (3.17) in \cite{TS} with zero
external field.\footnote{In their second equation, 
the range $1 \le i < \alpha$ should be
corrected as $1 \le i \le \alpha$. Also in their third equation
$d_1$ should be replaced with $d_i$.}

To obtain the free energy per site recall that 
$T^{(1)}_1(u,v)$ satisfies 
the inversion identity
$$
T^{(1)}_1(v+i) T^{(1)}_1(v-i) =T_0(v+2i) T_0(v-2i) (1+Y^{(1)}_1(v)).
$$
Again, the ANZC property of the both sides leads to
\begin{eqnarray*}
\ln T^{(1)}_1(u,v) &=&
    s_1*\ln(1+Y^{(1)}_1)(v)  +
     \ln \phi(v+i(u+2)) \phi(v-i(u+2))    \\
    &+& N \int^\infty_{-\infty} \frac{dk}{2k}e^{ikv}
       \frac{\sh ku \sh(1-\frac{\pi}{\theta})k} 
                  {\ch k \sh\frac{\pi k}{\theta}}.
\end{eqnarray*}
Calculating the limit in (\ref{eqn:ebaxter})
we obtain
$$
f =-\frac{2\pi J\sin\theta}{\theta}\int^\infty_{-\infty} a_1(v)s_1(v) dv -
   k_B T \int^\infty_{-\infty} s_1(v) \ln(1+\eta^{(1)}_1(v)) dv, 
$$
where
\[
 a_1(v)=\frac{1}{2p_0}\frac{\sin\theta}{\ch(\theta v)-\cos\theta}.
\]
This coincides with eq.(3.12) in \cite{TS} under the convention
$k_B = 1$.

\section{Correlation length}
Let us study the correlation lengths of 
$\langle\sigma_j^{+}\sigma_i^{-}\rangle$ 
and
$\langle\sigma_j^{z}\sigma_i^{z}\rangle$ 
along the scheme (\ref{cl}).
They are relevant to the second and the third largest 
eigenvalues $T^{(2)}_1(u,v)$ and $T^{(3)}_1(u,v)$ of the QTM,
respectively.
The former lies in the sector 
$m=N/2-1$ and the latter in $m=N/2$, where $m$ is the number 
of the Bethe ansatz roots in (\ref{qdefinition}). 
In this section we shall exclusively consider the case 
$p_0 \in {\bf Z}_{\ge 3}$ and $J > 0$, when 
the $Y$-system and $Y$-functions take the simple forms 
(\ref{ysystemc1})--(\ref{yfnc3}).

First we need to allocate the zeros of 
$T_{n-1}^{(k)}(u,v)$ ($k=2,3$) for $2 \le n \le p_0$
in the complex $v$-plane when $u$ is negative small.
Based on numerical studies, we have the following 
for $u$ negative small (typically $u \sim -0.1$).
\begin{cj}
For $2 \le n < p_0$, $T^{(2)}_{n-1}(u,v)$ has two real zeros 
$\pm \zeta^{(2)}_{n-1}$ for some 
$\zeta^{(2)}_{n-1} \in {\bf R}_{> 0}$. 
All the other zeros of
$T^{(2)}_{n-1}(u,v)$ ($2 \le n \le p_0$) are located 
on the almost straight lines 
$\Im v =  \pm n \hbox{ mod } 2p_0$. 
\end{cj}
For example see Fig.\ref{zeros23} showing the zeros of 
$T^{(2)}_{n-1}(u,v)$ for the case $p_0=5$, $u=-0.1$, $n=2,3,4,5$ and 
$N=20$.
The main difference from the largest eigenvalue case is 
the presence of the two real zeros for $n < p_0 $.
Their absence for $n=p_0$ can be explained as follows.
$T_{n-1}^{(2)}(u,v)$ in (\ref{dvf})--(\ref{qdefinition}) is a 
Laurent polynomial of $e^{\frac{\theta}{2}v}$.
When $m = N/2-1$, its highest/lowest terms are 
proportional to 
$\frac{\sin n\theta}{\sin \theta}e^{\pm \frac{N\theta}{2}v}$.
This is vanishing when $n = p_0$, therefore the 
number of zeros decreases from $N$ to $N-2$.
As a result  $Y^{(2)}_{p_0-2}(v)$ tends to zero 
as $e^{-\theta |v|}$ for $v \rightarrow \pm\infty$.

As for the third largest eigenvalue we have the following 
for $u$ negative small (typically $u \sim -0.1$).
\begin{cj}
$T^{(3)}_{n-1}(u,v)$ has two real zeros 
$\pm \zeta^{(3)}_{n-1}$ for some 
$\zeta^{(3)}_{n-1} \in {\bf R}_{> 0}$ for $n < p_0$ and 
a double zero at $\zeta^{(3)}_{p_0-1} = 0$ for $n = p_0$.
All the other zeros of
$T^{(3)}_{n-1}(u,v)$ ($2 \le n \le p_0$) are located 
on the almost straight lines
$\Im v =  \pm n \hbox{ mod } 2p_0$. 
\end{cj}
See Fig.\ref{zeros23} showing the zeros of 
$T^{(3)}_{n-1}(u,v)$ under the same conditions with
$T^{(2)}_{n-1}(u,v)$.
Again the main difference from the largest eigenvalue is the 
two additional zeros on the real axis.

When $v\to \pm\infty$, the $Y$-functions 
$Y^{(k)}_j$ and $K^{(k)}$ built from $T^{(k)}_{n-1}$ via 
(\ref{yfnc1})--(\ref{yfnc3}) have the asymptotic values
\bea
 Y^{(2)}_j &\rightarrow& 
         \frac{\sin\frac{\pi(j+2)}{p_0}\sin\frac{\pi j}{p_0}}
         {\sin^2\frac{\pi}{p_0}},\quad
 Y^{(3)}_j \rightarrow j(j+2), \label{yasympt} \\
 K^{(2)}&\rightarrow&-1,\quad \qquad \qquad \qquad K^{(3)}\rightarrow
 p_0-1.\label{kasympt}
\eea
To apply Lemma 1 to
the $Y$-system
(\ref{ysystemc1})-(\ref{ysystemc3}), 
we modify the $Y$-functions as
\bea
 \widetilde{Y}^{(k)}_j(v)&=&\frac{Y_j^{(k)}(v)}{F_j^{(k)}(v)} 
\quad \mbox{for}\,\, 1\le j \le p_0-2,\\
 \widetilde{K}^{(k)}(v)&=& \frac{K^{(k)}(v)}
 {\tnh\frac{\pi}{4}(v+\zeta^{(k)}_{p_0-2})
 \tnh\frac{\pi}{4}(v-\zeta^{(k)}_{p_0-2})},
\eea
where
\bean
F_1^{(k)}(v)&=&\{\tnh\frac{\pi}{4}(v+i(1+u))
           \tnh\frac{\pi}{4}(v-i(1+u))\}^{\frac{N}{2}}
           (\tnh\frac{\pi}{4}v)^{2\delta_{k,3}}g^{(k)}(v)
\quad \mbox{for} \,\,p_0=3, \\
F_1^{(k)}(v)&=&\tnh\frac{\pi}{4}(v+\zeta^{(k)}_{2})
           \tnh\frac{\pi}{4}(v-\zeta^{(k)}_{2})
          \{\tnh\frac{\pi}{4}(v+i(1+u))
           \tnh\frac{\pi}{4}(v-i(1+u))\}^{\frac{N}{2}}
\quad \mbox{for} \,\, p_0\ne 3,\\
F_j^{(k)}(v)&=&
          \tnh\frac{\pi}{4}(v+\zeta^{(k)}_{j+1})
           \tnh\frac{\pi}{4}(v-\zeta^{(k)}_{j+1})
          \tnh\frac{\pi}{4}(v+\zeta^{(k)}_{j-1})
           \tnh\frac{\pi}{4}(v-\zeta^{(k)}_{j-1})
\quad\mbox{for} \,\, 2\le j \le p_0-3, \\
F_{p_0-2}^{(k)}(v)&=&\tnh\frac{\pi}{4}(v+\zeta^{(k)}_{p_0-3})
           \tnh\frac{\pi}{4}(v-\zeta^{(k)}_{p_0-3})
          (\tnh\frac{\pi}{4}v)^{2\delta_{k,3}}g^{(k)}(v)
\quad\mbox{for} \,\, p_0 \ne 3.
\eean
The factor $g^{(k)}(v)$ defined by 
\[
g^{(k)}(v)=\cases{
             \exp\left(-\frac{\pi v}{p_0}\tnh\frac{\pi}{4}v\right) 
              & for $k=2$ \cr
              1 & for $k=3$ \cr}
\]
has been included to compensate the singularity caused by
$Y_{p_0-2}(v)$ tending to zero as $e^{-\frac{\pi}{p_0}|v|}$ at
$v=\pm\infty$.\footnote{This is a distinct feature of the present case
compared with \cite{Fen,BLZ,R}.}
The zeros $v = \zeta^{(k)}_j$ 
of $T^{(k)}_j(u_N,v)$ depend on $N$ and 
converge to some finite values in the Trotter limit $N\to\infty$.
By abuse of notation we shall also write
their limit as $\zeta^{(k)}_j$.
($\zeta^{(3)}_{p_0-1}=0$ is valid irrespective of $N$.)
Proceeding as in the free energy case, we get the 
non-linear integral equations obeyed by $\eta^{(k)}_j$ and $\kappa^{(k)}$:
\bea
\ln \eta^{(k)}_1(v)&=&
-\frac{\beta\pi J \sin\theta}{2\theta\ch(\frac{\pi v}{2})}+
  s_1\ast\ln\left((1+\kappa^{(k)})^2 h^{(k)}\right)(v)
+ \delta_{k,2}\left(\pi i-\frac{\pi}{3}v\tnh\frac{\pi v}{4}
  \right) \nn \\ 
&&+2\delta_{k,3}\ln\tnh\frac{\pi}{4}v
\quad \mbox{for}\,\,p_0=3,\label{nle1}\\
\ln \eta_1^{(k)}(v)&=&
-\frac{\beta\pi J \sin\theta}{2\theta\ch(\frac{\pi v}{2})}+
    s_1\ast\ln(1+\eta^{(k)}_2)(v) 
+ \ln \left\{ \tnh\frac{\pi}{4}(v+\zeta^{(k)}_2)
  \tnh\frac{\pi}{4}(v-\zeta^{(k)}_2) \right\}\quad \mbox{for}\,\, 
p_0\ne 3,\nn \\
&& \\ 
\ln \eta^{(k)}_j(v)&=&s_1\ast\ln(1+\eta^{(k)}_{j-1})(1+\eta^{(k)}_{j+1})(v) 
+  \ln\left\{\tnh\frac{\pi}{4}(v+\zeta^{(k)}_{j+1})
        \tnh\frac{\pi}{4}(v-\zeta^{(k)}_{j+1})\right\}\nn \\
&&+
      \ln\left\{\tnh\frac{\pi}{4}(v+\zeta^{(k)}_{j-1})
        \tnh\frac{\pi}{4}(v-\zeta^{(k)}_{j-1})\right\}
\quad \mbox{for}\,\, 2\le j\le p_0-3,\\ 
\ln \eta^{(k)}_{p_0-2}(v)&=&s_1\ast\ln\left((1+\eta^{(k)}_{p_0-3})
(1+\kappa^{(k)})^2 h^{(k)}\right)(v)
+\ln\left\{\tnh\frac{\pi}{4}(v+\zeta^{(k)}_{p_0-3})
                     \tnh\frac{\pi}{4}(v-\zeta^{(k)}_{p_0-3})\right\}\nn \\
&& +2\delta_{k,3} \ln\tnh\frac{\pi}{4}v
+\delta_{k,2}\left(\pi i -\frac{\pi}{p_0}v\tnh\frac{\pi v}{4}\right)
\quad \mbox{for} \,\,p_0\ne 3, \label{extba}\\
\ln\kappa^{(k)}(v)&=&s_1\ast\ln(1+\eta^{(k)}_{p_0-2})(v)
              + \ln \left\{
     \tnh\frac{\pi}{4}(v+\zeta^{(k)}_{p_0-2})
                    \tnh\frac{\pi}{4}(v-\zeta^{(k)}_{p_0-2})
                    \right\}+\delta_{k,2}\pi i, \label{nle2}
\eea
where 
\[
h^{(k)}(v)=\cases{
             \exp\left(\frac{2\pi}{p_0}\left(v \tnh\frac{\pi v}{2}-\frac{1}
             {\ch\frac{\pi v}{2}}\right)\right) 
              & for $k=2$ \cr
              1 & for $k=3$ \cr}.
\]
Here the integration constants have been fixed from the asymptotic
values (\ref{yasympt}) and (\ref{kasympt}).
In addition to these we need to impose the consistency condition
coming from $T^{(k)}_j(\zeta^{(k)}_j) = 0$, which determines
the real zeros 
$\{\pm\zeta^{(k)}_j\,|\,\zeta^{(k)}_j>0, k=2,3, j\in\{1,\cdots,p_0-2\}\}$.
($\zeta^{(3)}_{p_0-1}=0$.)
{}From (\ref{yfnc2}) and (\ref{etadef}) this leads to 
setting $\eta^{(k)}_j(\zeta_j^{(k)}\pm i)=-1$ 
in (\ref{nle1})--(\ref{extba}).
Explicitly they read
\bea
&\hbox{for}&p_0=3, \nn \\
&&i \frac{\beta\pi J \sin\theta}{2\theta\sh(\frac{\pi\zeta^{(k)}_1}{2})}+
   s_1\ast\ln\left((1+\kappa^{(k)})^2 h^{(k)}\right)(\zeta^{(k)}_1+i) 
 -  \delta_{k,2}\frac{\pi}{3}(\zeta^{(k)}_1+i)
            \tnh\frac{\pi (\zeta^{(k)}_1+i)}{4} \nn \\
&&+\delta_{k,3}
       \left\{\pi i +\ln \left(\frac{\sh(\frac{\pi}{2}\zeta^{(3)}_1)+i}
{\sh(\frac{\pi}{2}\zeta^{(3)}_1)-i}\right)
\right\}=0,\label{z1}\\
&\hbox{for}& p_0\ne3, \nn \\
&&i\frac{\beta\pi J \sin\theta}{2\theta\sh(\frac{\pi \zeta^{(k)}_1}{2})}+
s_1\ast\ln(1+\eta^{(k)}_2)(\zeta_1^{(k)}+i)
+\ln \left(\frac{\sh(\frac{\pi}{2}\zeta^{(k)}_1)+
    i\ch(\frac{\pi}{2}\zeta^{(k)}_2)}
{\sh(\frac{\pi}{2}\zeta^{(k)}_1)-i\ch(\frac{\pi}{2}\zeta^{(k)}_2)}\right)
+\pi i=0,\nn \\
&& \\
&\hbox{for}&2\le j\le p_0-3, \nn \\
&& s_1\ast\ln(1+\eta^{(k)}_{j-1})(1+\eta^{(k)}_{j+1})(\zeta^{(k)}_j+i) \nn \\
&& + \ln \left(\frac{\sh(\frac{\pi}{2}\zeta^{(k)}_j)+
i\ch(\frac{\pi}{2}\zeta^{(k)}_{j+1})}
{\sh(\frac{\pi}{2}\zeta^{(k)}_j)-i\ch(\frac{\pi}{2}\zeta^{(k)}_{j+1}) }\right)+
\ln \left(\frac{\sh(\frac{\pi}{2}\zeta^{(k)}_j)+
i\ch(\frac{\pi}{2}\zeta^{(k)}_{j-1})}
{\sh(\frac{\pi}{2}\zeta^{(k)}_j)-i\ch(\frac{\pi}{2}\zeta^{(k)}_{j-1}) 
}\right)+\pi i =0,\\
&\hbox{for}&p_0\ne3, \nn \\
&& s_1\ast\ln\left((1+\eta^{(k)}_{p_0-3})(1+\kappa^{(k)})^2h^{(k)}
\right)(\zeta^{(k)}_{p_0-2}+i)
+\ln \left(\frac{\sh(\frac{\pi}{2}\zeta^{(k)}_{p_0-2})+
                i\ch(\frac{\pi}{2}\zeta^{(k)}_{p_0-3})}
                {\sh(\frac{\pi}{2}\zeta^{(k)}_{p_0-2})
               -i\ch(\frac{\pi}{2}\zeta^{(k)}_{p_0-3})}\right) \nn \\
&&+\delta_{k,3}
\left\{\pi i + \ln \left(\frac{\sh(\frac{\pi}{2}\zeta^{(k)}_{p_0-2})+
                              i}
                              {\sh(\frac{\pi}{2}\zeta^{(k)}_{p_0-2}) 
                              -i}\right)\right\}
- \delta_{k,2}\frac{\pi}{p_0}(\zeta^{(2)}_{p_0-2}+i)
                     \tnh\frac{\pi (\zeta^{(2)}_{p_0-2}+i)}{4}=0,\nn \\
&& \label{z2}
\eea
where the convolutions should be interpreted as
\[
 s_1*g(\zeta+i)=\mbox{p.v.}\left(\int^{\infty}_{-\infty}\frac{g(x)}
{4i\sh\frac{\pi}{2}(\zeta-x)}dx\right)+\frac{1}{2}g(\zeta).
\]
Here p.v. means the principal value.
Since $T^{(k)}_1(u_N,0)$ is negative from the numerical experiment, 
$\lim_{N\to\infty}\vert T^{(k)}_1(u_N,0) \vert$ can be expressed as
\bea
 \lim_{N\to\infty}\vert T^{(k)}_1(u_N,0)\vert&=& -\frac{J}{2}\beta\cos\theta+
                        \int^{\infty}_{-\infty}dv s_1(v)\ln(1+\eta^{(k)}_1(v))
                                      \nonumber \\
                                   &&+ \frac{2\pi J\sin\theta}{\theta}\beta
                                      \int^{\infty}_{-\infty}dv a_1(v)s_1(v)+ 
                                      2\ln\tnh\frac{\pi}{4}\zeta^{(k)}_1.
\eea
Finally we obtain the correlation length (\ref{cl}) as 
\be
\frac{1}{\xi_k}=-2\ln\tnh\frac{\pi}{4}\zeta^{(k)}_1-
               \int^{\infty}_{-\infty}dv s_1(v)\ln\left
               (\frac{1+\eta^{(k)}_1(v)}{1+\eta^{(1 )}_1(v)}\right).
\ee

We can solve (\ref{nle1})--(\ref{z2}) numerically as follows.
First we solve the BAE (\ref{bae}) numerically for a finite $N$
and determine the 
$Y$-functions and their real zeros. 
This serves as the first approximation of their large $N$ limit
$\eta^{(k)}_j$ and $\kappa^{(k)}$.
Second we input them 
into the rhs of (\ref{nle1})-(\ref{nle2}) and 
get the new $\eta$-functions in the lhs as an output.
Third we substitute the output $\eta$-functions 
into (\ref{z1})-(\ref{z2}). Solving them by Newton's method 
a new output for the zeros $\zeta^{(k)}_j$ can also be constructed.
Finally by iterating the second and the third processes
in the above until adequate
convergence is achieved, 
the $\eta$-functions and their real zeros are determined accurately.
In this way the present approach enables us to overcome
the difficulty of the naive numerical extrapolation 
of $T^{(k)}_1$ as $N \rightarrow \infty$ mentioned in the introduction.

For a comparison we depict the functions 
$\eta^{(k)}_2(v) \,(k=1,2,3)$ 
in Fig.\ref{etag23} under the same parameters.
We also include the graphs of the correlation lengths $\xi_2$ and $\xi_3$ 
for $p_0=2,3,4,5$ in Fig.\ref{correlation}.
In the low temperature limit these results agree with 
the earlier ones in  \cite{KZeit,KBIbook}.
\bea
 \lim_{\beta\to\infty}
\xi_2(\beta)/\beta&=&\frac{J\pi \sin \theta}{2(\pi-\theta)\theta},
\label{known1}\\
 \lim_{\beta\to\infty} 
\xi_3(\beta)/\beta&=&\frac{J(\pi-\theta)\sin\theta}{2\pi\theta}
\label{known2}
\eea
with high accuracy.
\section{ Summary and discussion}
We have revisited the thermodynamics of the spin
$\frac{1}{2}$ XXZ model at roots of unity 
by the QTM method.
Functional relations indexed by the TS numbers are 
found among the fusion hierarchy of QTM's ($T$-system) and 
their certain ratios ($Y$-system).
As a peculiar feature of a general root of unity,
the $Y$-functions (\ref{ydef1})--(\ref{ydef3}) and the 
$Y$-system (\ref{ysystem1})--(\ref{ysystem4}) are 
considerably involved compared with those in \cite{KNS1}.
Nevertheless they have a nice analyticity 
allowing a transformation to integral equations.
Our approach simplifies the numerics to examine the 
analyticity drastically in that
only the largest eigenvalue sector of the QTM $T_1$ is needed
for the free energy.
We have set up Conjecture 1 on the zeros of QTM's
supported by an extensive numerical study.
The resulting integral equations exactly 
coincide with the TBA equation in \cite{TS} based on the string hypothesis.
Another and more significant advantage of the present method is
to allow us to study correlation lengths on an equal footing 
with the free energy by considering other eigensectors of $T_1$.
The additional zeros and poles coming into the ANZC strips play 
a fundamental role in characterizing the relevant excited states.
We have considered the second and the third largest 
eigenvalues of $T_1$, which are related to the 
spin--spin correlation lengths for 
$\langle\sigma_j^{+}\sigma_i^{-}\rangle$ and 
$\langle\sigma_j^z\sigma_i^z\rangle$, respectively.
The excited state TBA equation is derived and numerically solved 
to evaluate the correlation lengths.
The result shows a good agreement with the earlier one 
in the low temperature limit.

Let us remark a few straightforward generalizations
of the present results.
(1) the XYZ model, (2) higher spin cases
and (3) inclusion of external field $h$.
For (1) and (2), the $T$ and $Y$-systems remain 
essentially the same. 
We have an additional periodicity in the real direction 
in the complex $v$-plane for  (1).
This does not complicate actual calculations too much. 
In  (2) the driving term will enter the TBA equation in a different manner from 
(\ref{nlie1})-(\ref{nlie3}).
As noted in \cite{KR1}, the commensurability 
between the magnitude of the spin and the anisotropy
parameter would be of issue.
This is also an interesting problem in view of the present approach.
For  case (3) 
the BAE (\ref{bae}) should be modified 
with an extra $\exp(2\beta h)$ factor in the rhs.
Consequently, the BAE roots for the largest eigenvalue
will distribute away from the real axis.
This is a significant difference from the usual row-to-row 
case where they remain on the real axis even for $h \neq 0$.
The numerical check of the ANZC property therefore
needs more elaboration.
The $T$-system (\ref{tsystem2}) also needs to be modified into
\[
 T_{y_{\alpha+y_{\alpha-1}-1}}(v)=T_{y_\alpha-y_{\alpha-1}-1}(v)
 +2(-1)^{mz_\alpha}\ch(\beta h y_\alpha)T_{y_{\alpha-1}-1}(v+iy_\alpha).
\]
Correspondingly, (\ref{ysystem3}) is replaced by 
\[
 Y_{m_\alpha-1}(v)=K(v)^2+2\ch(\beta h y_\alpha)K(v).
\] 
These modifications are consistent with \cite{TS} from the 
string hypothesis. Explicit evaluation of the effects of the magnetic field
on correlation lengths will be an interesting problem
manageable within the present scheme.
\vskip0.5cm\noindent
{\bf Acknowledgements}.
A.~K. thanks M.~T.~Batchelor, E.~A. and R.~J.~Baxter and 
V.~V.~Bazhanov for hospitality at {\em International workshop on
statistical mechanics and integrable systems},
July 20 -- August 8, 1997 held in Coolangatta and Canberra,
where part of this work was presented.
He also thanks A.~Berkovich, B.~M.~McCoy and 
A.~Schilling for comments.
A.~K. and K.~S. are grateful to M.~Takahashi for useful discussions.

\appendix
\rnc{\theequation}{A.\arabic{equation}}\setcounter{equation}{0}
\section*{Appendix A \,\,Takahashi-Suzuki (TS) numbers} \label{app:TSnumbers}
Given $\{\nu_j\}$ in the continued fraction expansion 
(\ref{p0exp}) we define the sequences of numbers
$\{m_j\}_{j=0}^{\alpha+1},
\{p_j\}_{j=0}^{\alpha+1}, 
\{y_j\}_{j=-1}^{\alpha}, 
\{z_j\}_{j=-1}^{\alpha}, \{n_j\}_{j\ge 1}, 
\{\widetilde{ n}_j\}_{j\ge 1}$ and $\{w_j\}_{j\ge 1}$ as follows.
The sequence $\{m_j\}_{j=0}^{\alpha+1}$ is define by
\begin{eqnarray}
&&m_j = \nu_1 + \nu_2 + \cdots + \nu_j \quad 0 \le j \le \alpha,
\label{msec}\\
&&m_{\alpha+1} = \infty. \nonumber
\end{eqnarray}
The sequence $\{ p_j \}_{j=0}^{\alpha+1}$ is defined by
\begin{eqnarray}\label{pdefinition}
&&p_j = p_{j-2} - \nu_{j-1} p_{j-1} \quad 2 \le j \le \alpha+1,\nonumber\\
&&p_0 = \frac{\pi}{\theta},\, p_1 = 1, p_2 = p_0 - \nu_1.\nonumber
\end{eqnarray}
It can be easily shown that
\begin{eqnarray}
&&p_{\alpha+1} = 0,\label{eq:p1}\\
&&p_j < \frac{p_{j-1}}{\nu_j},\, p_j < \frac{p_0}{2}\quad
\mbox{if } 1 \le j \le \alpha+1,\label{eq:p2}\\
&&2p_j + 2p_{j+1} < p_0\quad \mbox{if } 2 \le j \le \alpha \hbox{ or }
j=1, \nu_1 \ge 3.\label{eq:p3}
\end{eqnarray}
The sequences $\{y_j\}_{j=-1}^\alpha$ and 
$\{z_j\}_{j=-1}^\alpha$ are defined by
\begin{eqnarray}
&&y_j = y_{j-2} + \nu_j y_{j-1} \quad 1 \le j \le \alpha,\nonumber\\
&&y_{-1} = 0, y_0 = 1, y_1 = \nu_1, y_2 = 1 + \nu_1\nu_2,
\label{yzdefinition1}\\
&&z_j = z_{j-2} + \nu_j z_{j-1} \quad 1 \le j \le \alpha,\nonumber\\
&&z_{-1} = 1, z_0 = 0, z_1 = 1, z_2 = \nu_2.\label{yzdefinition2}
\end{eqnarray}
Obviously, $z_j = y_{j-1}\vert_{\nu_k \rightarrow \nu_{k+1}}$ and 
they are all positive integers except $y_{-1} = z_0 = 0$.
By induction one can verify 
\begin{eqnarray}
&&y_j = z_j p_0 + (-1)^jp_{j+1} \quad  -1 \le j \le \alpha, 
\label{yzproperty1}\\
&&y_\alpha = z_\alpha p_0,
\label{yzproperty2}
\end{eqnarray}
where the latter is a consequence of the former with $j=\alpha$ and 
(\ref{eq:p1}).
In fact $GCD(y_\alpha, z_\alpha) = 1$ is valid.
Now we introduce the Takahashi-Suzuki (TS) numbers
$\{ n_j \}_{j\ge 1}$ \cite{TS} and their slight
rearrangement
$\{ {\widetilde n_j}\}_{j \ge 1}$ by
\begin{eqnarray}
&&n_j = y_{r-1} + (j-m_r)y_r \quad m_r \le j < m_{r+1}, \label{def:TSn}\\
&&{\widetilde n}_j = y_{r-1} + (j-m_r)y_r \quad m_r < j \le m_{r+1}. 
\label{ndefinition}
\end{eqnarray}
Obviously, ${\widetilde n}_j = n_j$ except 
${\widetilde n}_{m_r} = y_r$ while $n_{m_r} = y_{r-1}$.
In particular, there is a duplication $n_1 = n_{m_1} = 1$, while
the modified sequence ${\widetilde n}_j$ is strictly increasing with $j$.
In this paper we are concerned with the first $m_\alpha+1$ of them.
As the set with multiplicity
$$
\{ {\widetilde n}_j \}_{j=1}^{m_\alpha + 1}  = 
\{ n_j \}_{j=1}^{m_\alpha + 1} \sqcup \{ y_\alpha \} \setminus \{1\}.
$$
We note that if $p_0 \ge 2$ we always have
${\tilde n}_j = j$ for $j = 1, 2$ and 3.
In parallel with (\ref{yzdefinition1}) and 
(\ref{yzdefinition2}) we consider the ``$z$-analogue''
$\{ w_j \}_{j \ge 1}$ of $\{n_j \}_{j\ge 1}$:
\begin{equation}\label{wdefinition}
w_j = z_{r-1} + (j - m_r)z_r - 1\quad m_r \le j < m_{r+1}. \nonumber
\end{equation}
For example, $w_1 = w_{m_2} = 0$ and $w_{m_1} = -1$.
It is possible to show 
\begin{equation}\label{wproperty}
\left[ \frac{n_j-1}{p_0} \right] = w_j + \delta_{j m_1}
\quad 1 \le j \le m_\alpha,
\end{equation}
where $[ x ]$ denotes the largest integer not exceeding $x$.
As a result the sequence $\{w_j \}$ is related to
the parity $v_j$ in (2.14) of \cite{TS} by
$v_j = (-1)^{w_j}$ for all $1 \le j \le m_\alpha$.
Using (\ref{yzproperty1}), (\ref{ndefinition}) and 
(\ref{wdefinition}) one can show 
\begin{equation}
\{w_jp_0 \pm {\tilde n}_{j+1} \pm y_r \hbox{ mod } 2p_0 \} 
= \{ p_0 \pm (p_r - (j+1 \pm 1 - m_r)p_{r+1}) \hbox{ mod } 2p_0 \}
\label{tsukatte}
\end{equation}
for $m_r \le j < m_{r+1}$. Here the signs $\pm$ are independent.

It is well known \cite{TS,KR1,Kore} that for the TS
numbers $n_j$, the equivalent conditions
\begin{eqnarray}
&&(-1)^{w_j} \sin\left(\frac{\pi k}{p_0}\right)
           \sin\left(\frac{\pi(n_j - k)}{p_0}\right) > 0,
\label{TScondition1}\\
&&\left[ \frac{k}{p_0} \right] + \left[ \frac{n_j-k}{p_0} \right]
= \left[ \frac{n_j-1}{p_0} \right],
\label{TScondition2}
\end{eqnarray}
hold for $k = 1, 2, \ldots, n_j - 1$.
It is interesting to observe the condition 
(\ref{TScondition1}) in the light of the 
associated fusion transfer matrix 
$T_{n_j-1}(u,v)$ (\ref{ftdefinition}).
{}From (\ref{fbw}) and (\ref{boxv}), we see that 
(\ref{TScondition1}) ensures that $\epsilon_{j j'}$ 
and $\epsilon'_{j j'}$ can be independent of their
indices for the constituent fusion Boltzmann weights to be real.

\rnc{\theequation}{B.\arabic{equation}}\setcounter{equation}{0}
\section*{Appendix B \,\,ANZC property of $Y_j^{(1)}(v)$} 
Let us check the applicability of Lemma 1 in section 5
to the $Y$-system (\ref{ysystem1})--(\ref{ysystem4}) 
by admitting Conjecture 1.
Apart from the exceptional 
Case 1, 2 and 3 listed there,
we are to verify that $Y^{(1)}_{j}(v)$, for example, is 
ANZC in the strip ${\cal S}[x]$ whenever the 
combination $Y^{(1)}_{j}(v+ix)Y^{(1)}_{j}(v-ix)$ takes place.
Case 1 has been argued in section 5 and Case 2 and 3 will be considered in 
Appendix C.

Conjecture 1 tells that the zeros and poles of the 
$Y$-functions (\ref{ydef1})--(\ref{ydef3}) are located as
\bea
Y^{(1)}_j(v), 1+Y^{(1)}_j(v) :&& \Im v 
\sim w_jp_0\pm\widetilde{n}_{j+1}\pm y_r 
\quad 1 \le j \le m_\alpha -1,\nonumber\\
K^{(1)}(v):&&  \Im v 
\sim w_{m_\alpha-1}p_0 \pm y_{\alpha-1} + y_{\alpha}. \nonumber
\eea
Here the signs $\pm$ are independent and we have taken the 
periodicity under $v \rightarrow v + 2p_0i$ into account.
See (\ref{tperiod}).
{}From (\ref{tsukatte}) these functions are ANZC 
in the following strips:
\bea
\hbox{For}&& m_r\le j\le m_{r+1}-2\,\,(0 \le r \le \alpha-1),\nn \\ 
           &&Y^{(1)}_j(v),1+Y^{(1)}_j(v):\,\,
             \st[p_0-p_r+(j-m_r)p_{r+1}-\delta],\label{anzc1}\\
           && Y^{(1)}_{m_r-1}(v),1+Y^{(1)}_{m_r-1}(v):\,\,
             \st[p_0-p_r-p_{r+1}-\delta], \label{anzc2}\\
           && 1+Y^{(1)}_{m_{\alpha}-1}(v), 1+K^{(1)}(v),K^{(1)}(v):\,\,
             \st[p_0-p_{\alpha}-\delta]\label{anzc5}.
\eea
%
In the above $\delta$ denotes a small real number
$(\vert \delta \vert \sim 10^{-2})$ caused by the deviations 
of the actual zeros of $T^{(1)}_{n-1}(v)$ from the straight line 
specified in Conjecture 1.
They of course depend on the $Y$-functions but have been 
denoted by the same symbol for the sake of simplicity.
On the other hand, Lemma 1 is applicable to the $Y$-system 
(\ref{ysystem1})--(\ref{ysystem4}) if 
the $Y$-functions are ANZC in the strips:
\bea
    &&Y^{(1)}_j(v):\,\,\st[p_{r+1}]\,\,(m_r\le j\le m_{r+1}-2, \, 
                        0 \le r \le \alpha-1),\label{anzc7}\\
           && Y^{(1)}_{m_r-1}(v):\,\,\st[p_r+p_{r+1}] \,
              (1 \le r \le \alpha-1),\label{anzc9}\\
           &&1+Y^{(1)}_{m_r-1}(v):\,\,\st[p_r-p_{r+1}]\,
              (1 \le r \le \alpha-1),\label{anzc10}\\
           && 1+Y^{(1)}_{m_r-2}(v):\,\,\st[p_{r+1}]\,
              (1 \le r \le \alpha-1),\label{anzc11} \\
           && 1+Y^{(1)}_{m_r}(v):\,\,\st[p_r]\,
              (1 \le r \le \alpha-1),\label{anzc12}\\  
            && K^{(1)}(v) :\,\,\st[p_{\alpha}],\label{anzc13}\\
            && 1+K^{(1)}(v), \{1+Y^{(1)}_j(v)\}_{j=1}^{m_\alpha - 1} 
   \hbox{ other than (\ref{anzc10})--(\ref{anzc12})}:\,\,
            \st[0^+]\label{anzc14},
\eea
where $\st[0^+]$ means the vicinity along the real axis
which can be arbitrarily thin.
If $\nu_1=2$, the $r=1$ case of (\ref{anzc11}) is void.

Except for Cases 1, 2 and 3 in section 5, 
it is straightforward to verify that the strips in 
(\ref{anzc7})--(\ref{anzc14}) are narrower than those
in (\ref{anzc1})--(\ref{anzc5}) for 
the corresponding functions.
As an illustration we prove here that 
$\st[p_0 - p_r+(j-m_r)p_{r+1}-\delta]
\supset \st[p_{r+1}]$ for the strips in (\ref{anzc1}) and (\ref{anzc7}). 
The rest is a similar exercise.
We only have to show the inequality
\be
p_0 - p_r+(j-m_r)p_{r+1}-\delta>p_{r+1}\quad\hbox{for}\,\,
m_r\le j\le m_{r+1}-2.
\label{ieq}
\ee
Though this is incorrect for $j=1$ (hence $r = 0$), 
this case corresponds to Case 1, for which the difficulty has been 
cleared in section 5  by a modification 
of a $Y$-function.
Now suppose $j \neq 1$.
It is enough to check (\ref{ieq}) only for $j=m_r\,\,(r \ge 1)$.
Making use of the properties in (\ref{eq:p2}) one has
\be
p_0-p_{r+1}-p_r - \delta \ge p_0-p_2-p_1-\delta  = \nu_1-1 - \delta.
\ee
By noting that $\nu_1\ge 2$ and $\vert \delta  \vert \ll 1$,
the last quantity is non-negative, proving (\ref{ieq}).
\rnc{\theequation}{C.\arabic{equation}}\setcounter{equation}{0}
\section*{Appendix C\,\, ANZC property of 
$Y^{(1)}_j(v)$ for exceptional cases}

Let us show that Lemma 1 can still be applied to the $Y$-system in Case 2 and 3 
in section 5 after suitable recombinations of the $Y$-functions.
For a function whose logarithmic derivative can be 
Fourier transformed we use the notation
$$
{\cal F}[f](k) =
 \frac{1}{2\pi } \int_{-\infty}^{\infty} \frac{d}{dv} \ln f(v) e^{-ikv} dv. 
$$

We start with Case 2.
Explicitly it reads
\begin{eqnarray}
&& Y^{(1)}_2(v-i(4-p_0)) Y^{(1)}_2(v+i(4-p_0)) Y^{(1)}_2(v-i(p_0-2)) 
Y^{(1)}_2(v+i(p_0-2))= \nonumber\\
&&(1+Y^{(1)}_1(v+i(p_0-3))) (1+Y^{(1)}_1(v-i(p_0-3)))
(1+Y^{(1)}_2(v+i(4-p_0))) \nn \\
&&\times(1+Y^{(1)}_2(v-i(4-p_0)))(1+Y^{(1)}_3(v+i))(1+Y^{(1)}_3(v-i)),
\label{case2}
\end{eqnarray}
where
\be
1+Y^{(1)}_3(v)=\frac{T^{(1)}_3(v+i(p_0-3))T^{(1)}_3(v-i(p_0-3))}
           {T^{(1)}_2(v+i(4-p_0))T^{(1)}_2(v-i(4-p_0))},
\ee
and the other functions are given by (\ref{yfnc1}) and 
(\ref{yfnc2}).
{}From the Case 2 conditions on
$\nu_1$--$\nu_3$,  we have 
$p_0=4-\epsilon$  with $0<\epsilon\ll 1$.  
Thus the ANZC argument can not be 
applied to some factors in (\ref{case2}). 
For example the $Y^{(1)}_2(v)$-function in the lhs has zeros or poles along 
$\Im v \simeq 2$.
They are outside of $\st[4-p_0]$ but can be within 
$\st[p_0-2]$.
Similarly zeros of the $1+Y^{(1)}_1(v)$ in the rhs 
lie along $\Im v \simeq 1$ which is in the strip
$\st[p_0-3]$.
A recipe here is to consider the combination  
\begin{equation}
X(v) = \frac{Y^{(1)}_2(v+i) Y^{(1)}_2(v-i)}{1+Y^{(1)}_1(v)}=
  \frac{T^{(1)}_3(v+i) T^{(1)}_3(v-i)}{T_0(v+4 i) T_0(v-4 i)},
\label{defx}  
\end{equation}
which is ANZC in $\st[p_0-3]$.
With this, (\ref{case2}) can be rewritten as
\begin{eqnarray}
& &X(v+i(p_0-3)) X(v-i(p_0-3))  \nonumber \\
& &=
(1+Y^{(1)}_2(v+i(p_0-4)))(1+Y^{(1)}_2(v-i(p_0-4)))(1+Y^{(1)}_3(v+i)) \nn \\
&& \times (1+Y^{(1)}_3(v-i)).
\label{xfuncrel}  
\end{eqnarray}
At this stage, the lemma applies to both (\ref{defx}) and
(\ref{xfuncrel})
giving
\begin{eqnarray*}
{\cal F}[ X](k)  &=& 2 \ch k {\cal F}[Y^{(1)}_2](k) - 
{\cal F}[1+Y^{(1)}_1](k), \\
{\cal F}[ X](k)  &=&  
\frac{\ch(4-p_0)k}{ \ch(p_0-3)k} {\cal F}[1+Y^{(1)}_2](k)   \\
 &+& \frac{\ch k}{ \ch(p_0-3)k} {\cal F}[1+Y^{(1)}_3](k).
\end{eqnarray*}
Eliminating ${\cal F}[X](k)$ from these and doing
the inverse Fourier transformation, we get (\ref{nlie2}).

Next we consider Case 3.
\begin{eqnarray}
&& Y^{(1)}_1(v+i(1-p_2)) Y^{(1)}_1(v-i(1-p_2)) Y^{(1)}_1(v+i(1+p_2))  
Y^{(1)}_1(v-i(1+p_2)) \nonumber \\
&& = (1+Y^{(1)}_2(v+i))(1+Y^{(1)}_2(v-i))
(1+Y^{(1)}_1(v+i(3-p_0))) \nn \\
&&\times (1+Y^{(1)}_1(v-i(3-p_0))),
\label{case3}
\end{eqnarray}
where
\begin{eqnarray}
1+Y^{(1)}_2(v)&=&\frac{T^{(1)}_2(v+i(p_0-2))T^{(1)}_2(v-i(p_0-2))}
                {T_0(v+i(3-p_0))T_0(v-i(3-p_0))},
\end{eqnarray}
and the other functions are given by 
(\ref{yfnc1}) and (\ref{yfnc2}). 
Now $p_0=2+p_2$.
$Y^{(1)}_1(v)$ has  zeros (resp. poles) at 
$v=\pm(1+u)i$ (resp. $v=\pm(1-u)i$)  $\in \st[1+p_2]$,
which prevents the direct application of Lemma 1
in the last two factors in the lhs of (\ref{case3})
for $J >0$ (resp. $J<0$).
This can be remedied by introducing 
$\widetilde{Y}^{(1)}_1(v) = 
Y^{(1)}_1(v)/(\tnh(\frac{\pi}{4}\frac{v-i(1 \pm u)}{1+p_2})  
 \tnh(\frac{\pi}{4}\frac{v+i(1 \pm u)}{1+p_2}))^{\pm N/2}$
as in section 5.
In the rhs of (\ref{case3}), there are also some factors 
possessing zeros or poles and making Lemma 1 inapplicable.
However the new combinations
\begin{eqnarray*}
G_1(v) &=& (1+Y^{(1)}_2(v))(1+Y^{(1)}_1(v-i(p_0-2))) \\
       &=& \frac{T^{(1)}_2(v+i(p_0-2))}{T^{(1)}_1(v-i(3-p_0))}
           \frac{T^{(1)}_2(v-i(p_0-2))}{T_0(v+i(4-p_0))}
           \frac{T^{(1)}_1(v-i(p_0-1))}{T_0(v-ip_0)},\\
G_2(v) &=& (1+Y^{(1)}_2(v))(1+Y^{(1)}_1(v+i(p_0-2)))\\
       &=& \frac{T^{(1)}_2(v+i(p_0-2))}{T^{(1)}_1(v+i(3-p_0))}
           \frac{T^{(1)}_2(v-i(p_0-2))}{T_0(v-i(4-p_0))}
           \frac{T^{(1)}_1(v+i(p_0-1))}{T_0(v+ip_0)}
\end{eqnarray*}
are free of these spurious zeros and poles and ANZC in $\Im v\in[0,1]$ and
$[-1,0]$, respectively.
With their aid (\ref{case3}) can be written as
\begin{eqnarray*}
& & Y^{(1)}_1(v+i(1-p_2)) Y^{(1)}_1(v-i(1-p_2)) 
\widetilde{Y}^{(1)}_1(v+i(1+p_2))  
\widetilde{Y}^{(1)}_1(v-i(1+p_2))\\
&& = G_1(v+i) G_2(v-i).
\end{eqnarray*}
Solving these relations as in Case 2, we obtain the solution in Fourier space,
\begin{eqnarray*}
{\cal F}[ Y^{(1)}_1](k) 
 &=& \frac{i N \sh ku}{2 \ch k}  \\
 &+& \frac{ \ch(p_0-3)k}{2\ch p_2 k \ch k } {\cal F}[1+Y^{(1)}_1](k)  +
     \frac{1}{2  \ch p_2 k} {\cal F}[1+Y^{(1)}_2](k),  
\end{eqnarray*}
which can be transformed back to (\ref{nlie2}).
\rnc{\theequation}{D.\arabic{equation}}\setcounter{equation}{0}
\section*{Appendix D\,\,  Free fermion case}
Here we consider the free energy and the correlation
lengths for the free fermion case
$\Delta=0$ ($p_0=2,\alpha=1$) in (\ref{eqn:Hamiltonian}). 
In this case we have 
$\phi(v+4i) = (-1)^{\frac{N}{2}}\phi(v)$ and 
$Q(v+4i) = (-1)^mQ(v)$ from (\ref{phidefinition}) and 
(\ref{qdefinition}).
Thus (\ref{dvf}) simplifies to
\be
T_1(v)=\varrho(v)\frac{Q(v+2i)}{Q(v)},
\label{freeT}
\ee
where 
\[
\varrho(u,v)=\phi \bigl(v-i(u+2)\bigr)\phi(v+iu)+
             (-1)^m\phi\bigl(v+i(u+2)\bigr)
                \phi(v-iu). 
\]
One can directly show 
$T_1(v+i)T_1(v-i)=(-1)^m\varrho(v+i)\varrho(v-i)$.
This rhs is a known function, which is a distinct feature of 
the free fermion model.
We find it convenient to introduce
\bean
\widetilde{T}^{(k)}_1(u,v)&=& \frac{T_1^{(k)}(u,v)}{\phi(v+i(u+2))
                                  \phi(v-i(u+2))} \\
                          &=&\left(\frac{\phi(v+iu)}
                             {\phi(v+i(u+2))}+(-1)^{m}
                             \frac{\phi(v-iu)}{\phi(v-i(u+2))}\right)
                             \frac{Q(v+2i)}{Q(v)}.
\eean
It satisfies 
\be
\widetilde{T}_1^{(k)}(u,v+i)
\widetilde{T}_1^{(k)}(u,v-i)= (X(u,v)^{\frac{1}{2}}+
                                    (-1)^{\frac{N}{2}-m} 
                                    X(u,v)^{-\frac{1}{2}})^{2},
\label{eqn:rtfree}
\ee
where
\[
X(u,v)=\frac{\phi(v+i(u-1))\phi(v-i(u-1))}
            {\phi(v+i(u+1))\phi(v-i(u+1))}.
\]

First we consider the free energy characterized
by the largest eigenvalue $T_1^{(1)}(u,v)$.
It lies in the sector $m=N/2$.
Since the function $\widetilde{T}^{(1)}_1(u,v)$ is ANZC for 
$v \in \st[1]$, we have
\be
2\ch k \,{\cal F}[\widetilde{T}^{(1)}_1](k)=
2{\cal F}\,[ X^{\frac{1}{2}}+X^{-\frac{1}{2}}](k).
\label{eqn:ftr}
\ee
See Appendix C for the notation ${\cal F}$.
By the inverse Fourier transformation and the identity
\be
\int_{-\infty}^{\infty}\frac{e^{-ikv}}{2\ch k}dk
    =\frac{\pi}{2\ch\frac{\pi v}{2}} = 2\pi s_1(v),\nonumber
\ee
we get 
\be
\ln \widetilde{T}^{(1)}_1(u,v)
= 2\left[s_1\ast\ln (X^{\frac{1}{2}}+
                    X^{-\frac{1}{2}})\right](v). 
\ee
See (\ref{nlie4}). Using the relations
\[
\lim_{N\to\infty}X(u_N,v)=\exp\left(
                           \frac{J\beta}{\ch \frac{\pi v}{2}}\right),
\quad 
\lim_{N\to\infty}\widetilde{T}^{(k)}_1(u_N,0)=\lim_{N\to\infty}
                                               T^{(k)}_1(u_N,0),
\]
we obtain the free energy per site $f$ as
\bea 
f&=&-\frac{1}{\beta} \lim_{N \rightarrow \infty} 
      \ln T^{(1)}_1(u_N,0) \nonumber \\
 &=&-\frac{2}{\pi\beta} \int_0^{\frac{\pi}{2}}\ln\left(
2\ch (\frac{J\beta}{2} \cos\eta)\right) d\eta \label{fffe}
\eea
in agreement with \cite{Kat}.  

Next we consider the correlation length $\xi_2$ for
$\langle\sigma_j^{+}\sigma_i^{-}\rangle$ which
is related to the second largest eigenvalue
$T_1^{(2)}(v)$.
This lies in the sector $m=N/2-1$.
{}From a numerical check,
$T_1^{(2)}(u,v)$ is ANZC for $v \in \st[1]$. 
Therefore we can calculate it in the same
way as $T_1^{(1)}(u,v)$
The only difference is 
\be
\ln \widetilde{T}^{(2)}_1(v) = 2[s_1\ast \ln(X^{\frac{1}{2}}-
X^{-\frac{1}{2}})](v)
\ee
due to  (\ref{eqn:rtfree}) with $m=N/2-1$. 
Thus we have
\be
 \lim_{N\to \infty} \ln \vert T_1^{(2)}(u_N,0) \vert =
 \frac{2}{\pi} \int_0^{\frac{\pi}{2}}\ln\left(
2\sh (\frac{\vert J\beta\vert }{2} \cos\eta)\right) d\eta. \nonumber
\ee
Combining this with (\ref{fffe}) and (\ref{cl}) 
we find the $\xi_2$ for
$\langle\sigma_j^{+}\sigma_i^{-}\rangle$:
\be
\frac{1}{\xi_2}=-\frac{2}{\pi} \int_0^{\frac{\pi}{2}}\ln\left(
2\tnh (\frac{\vert J\beta \vert}{2} \cos\eta)\right) d\eta. 
\label{ffxi2}
\ee
This result coincides with those in \cite{InSuz2},
\cite{Kat}--\cite{LSM}.

Finally we consider the correlation length $\xi_3$ 
for $\langle\sigma_j^{z}\sigma_i^{z}\rangle$
characterized by the third largest
eigenvalue $T_1^{(3)}(u,v)$.
This lies in the sector $m=N/2$,
which is the same with $T_1^{(1)}(u,v)$.
All the BAE roots for $T_1^{(1)}(u,v)$ are real 
solutions of $\varrho(v)=0$.
The set of the BAE roots for $T_1^{(3)}(u,v)$ is the same with 
the one for $T_1^{(1)}(u,v)$ 
except that the largest magnitude ones $\pm \zeta$ are replaced by 
0 and $2i$.
It follows from (\ref{freeT}) that 
$T^{(3)}_1(u,v) = T^{(1)}_1(u,v) 
\tnh\frac{\pi}{4}(v+\zeta)\tnh\frac{\pi}{4}(v-\zeta)$.
In the Trotter limit, the real zeros $\pm \zeta$ are 
the largest magnitude solutions to $\varrho(u_N, v)|_{N\to\infty}=0$ 
given by
\[
 \zeta=\frac{2}{\pi}\sh^{-1}\left(\frac{J\beta}{\pi}\right).
\]  
Thus we have
$ \lim_{N\to \infty}\ln \left|T_1^{(3)}(u_N,0)\right|=
 \lim_{N\to\infty}\ln \left|T_1^{(1)}(u_N,0)\right|
+2\ln\tnh\frac{\pi\vert \zeta \vert}{4}$, 
and obtain the $\xi_3$ for  
$\langle\sigma_j^{z}\sigma_i^{z}\rangle$ as
\bea
\frac{1}{\xi_3}&=&-2\ln\tnh\left(\frac{1}{2}
        \sh^{-1}\left(\frac{\vert J\beta\vert}{\pi}\right)\right) 
\nonumber \\
     &=&2 \sh^{-1}\left(\frac{\pi}{\vert J\beta\vert}\right). 
\label{ffxi3}
\eea
This agrees with \cite{Kat}--\cite{LSM}.

\begin{figure}[p]
  \centerline{\epsfig{file=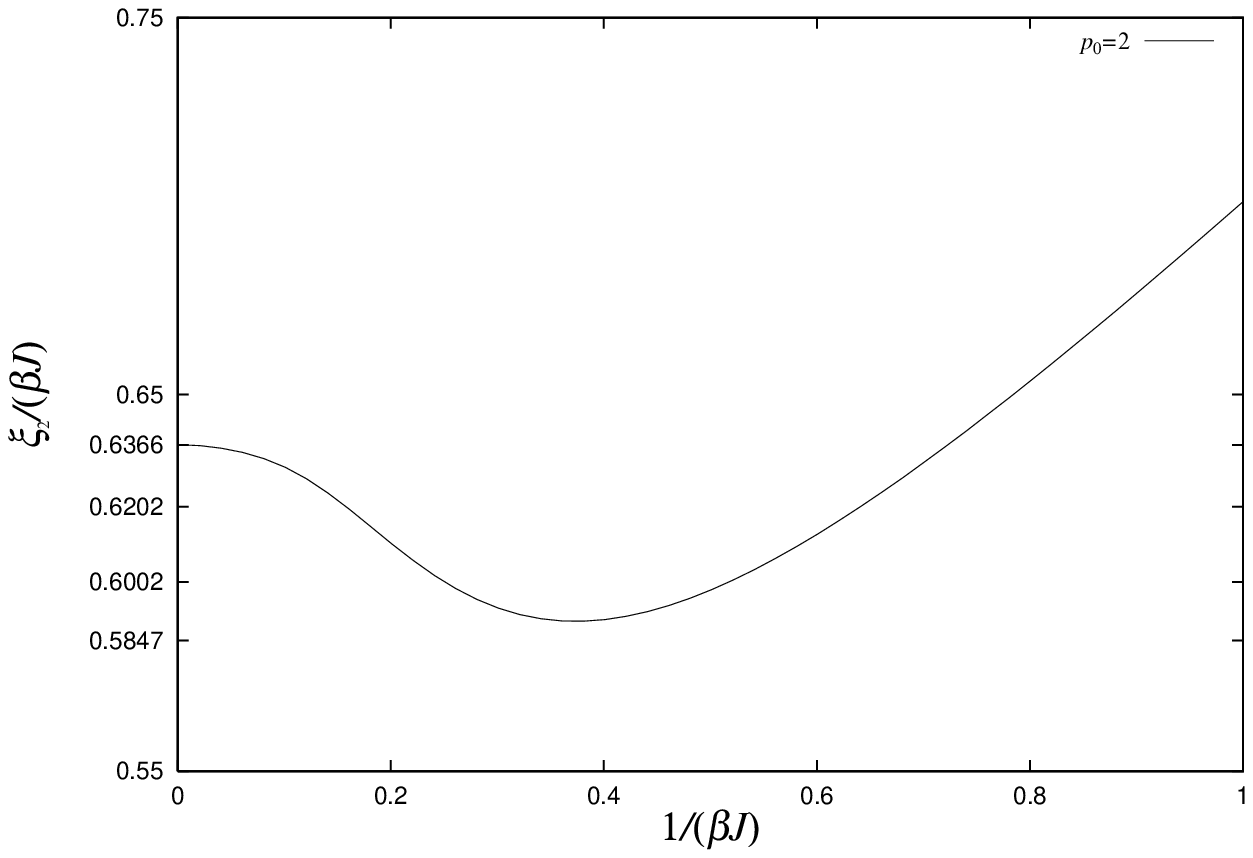}}   
  \centerline{\epsfig{file=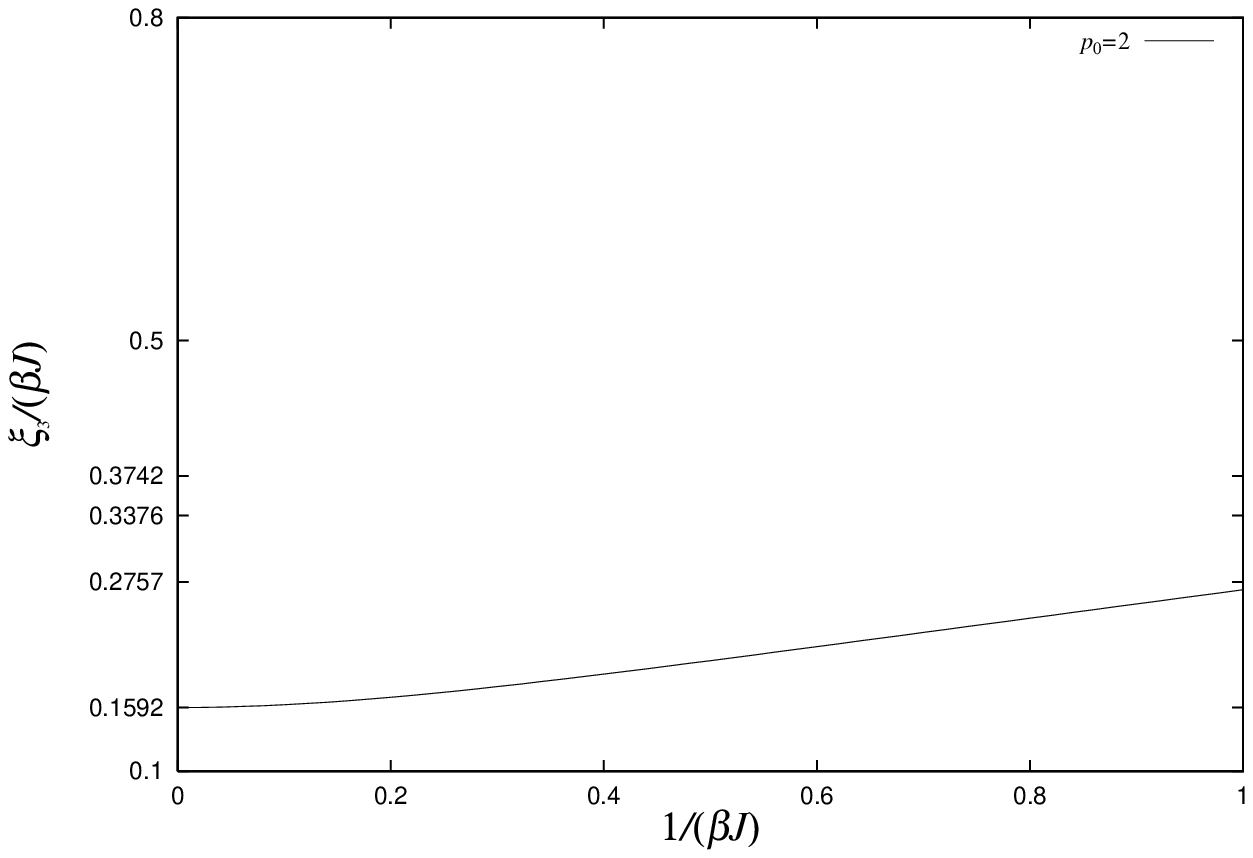}}
\caption{Ratio of the correlation length and the inverse 
temperature.
$\langle\sigma^+_j \sigma^-_i\rangle$ 
(upper) and 
$\langle\sigma^z_j \sigma^z_i\rangle$ (lower) for 
$p_0=3,4,5$ and the free fermion case $p_0=2$.
The known result (6.16)--(6.17) 
in the low temperature limit is also depicted by the symbol $\ast$.}
\label{correlation}
\end{figure}

The results (\ref{ffxi2}) and (\ref{ffxi3}) are also plotted 
in Fig.\ref{correlation}.

\begin{figure}[p]
  \centerline{\epsfig{file=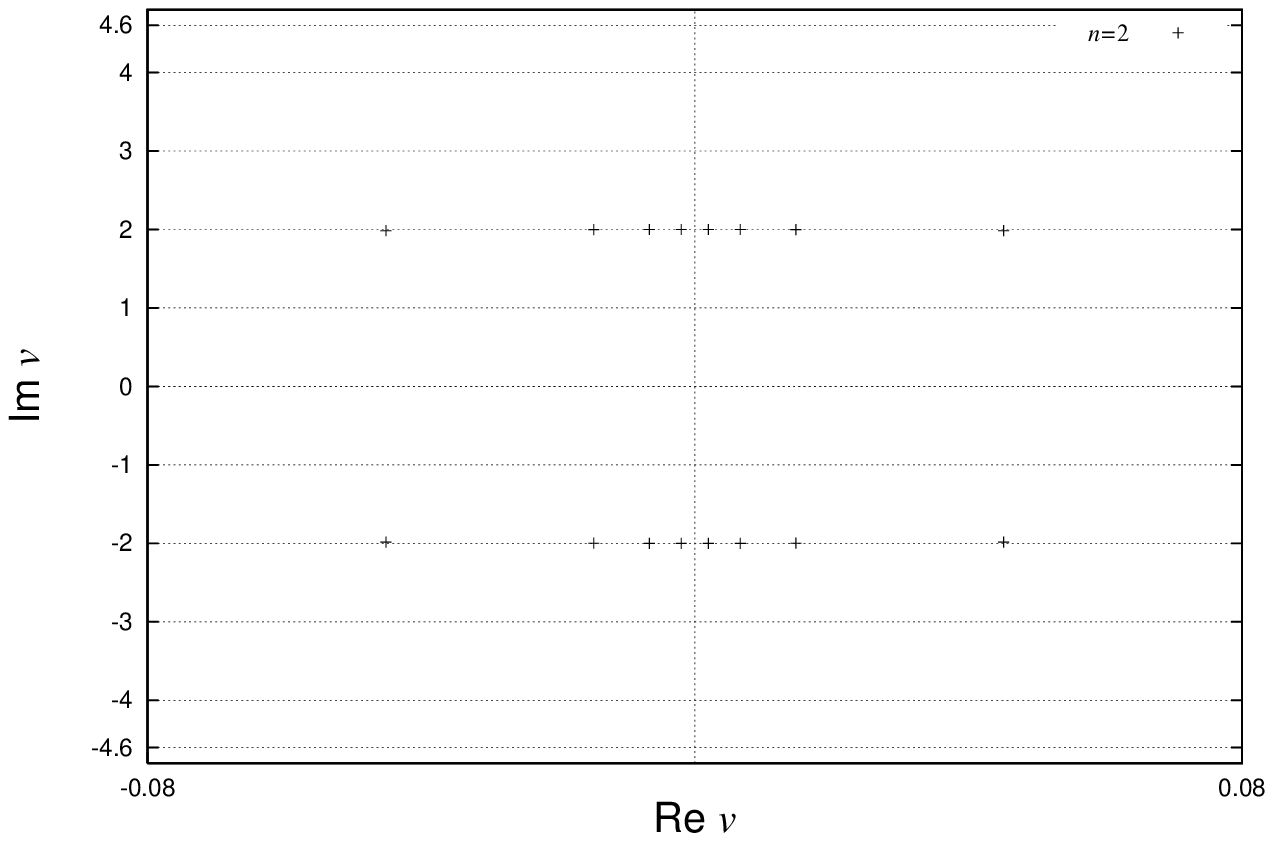,height=7cm}}
  \centerline{\epsfig{file=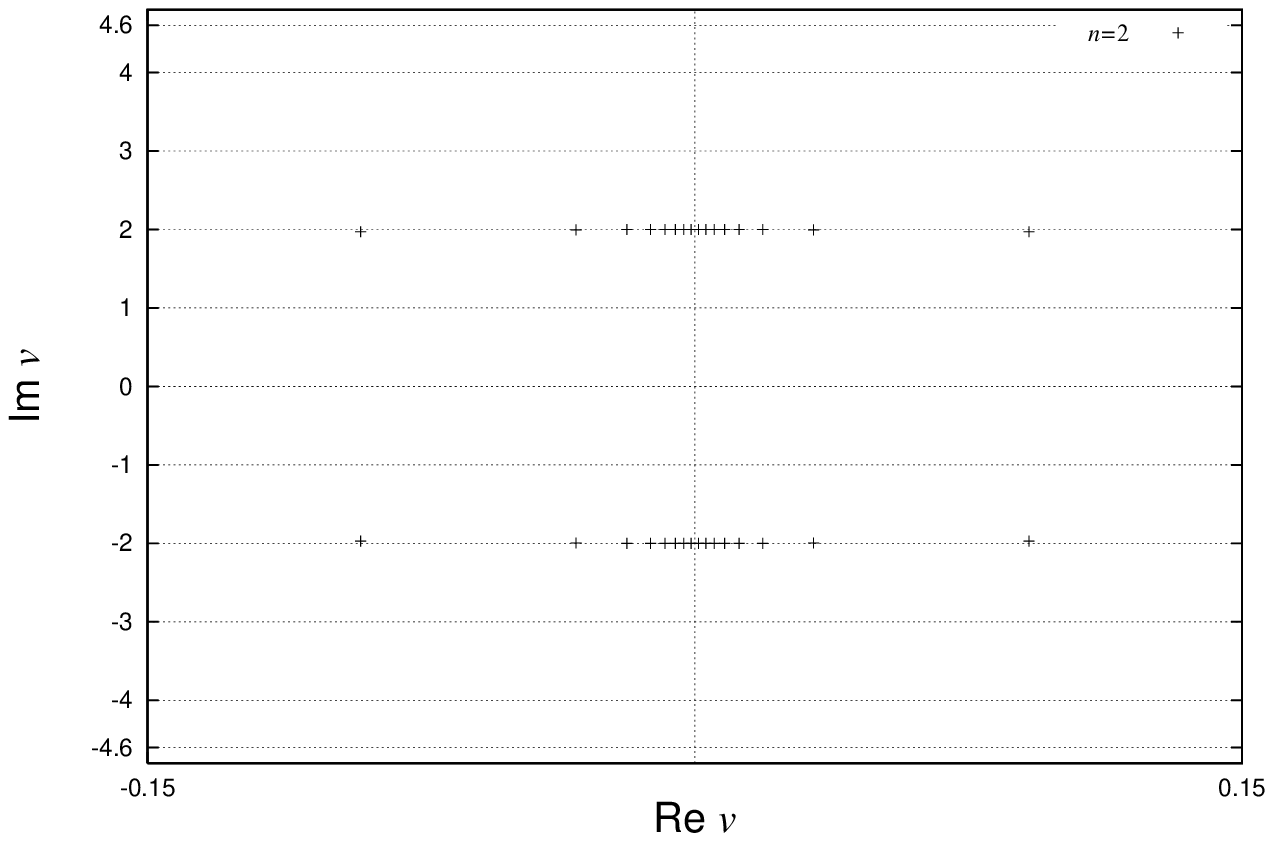,height=7cm}}
  \centerline{\epsfig{file=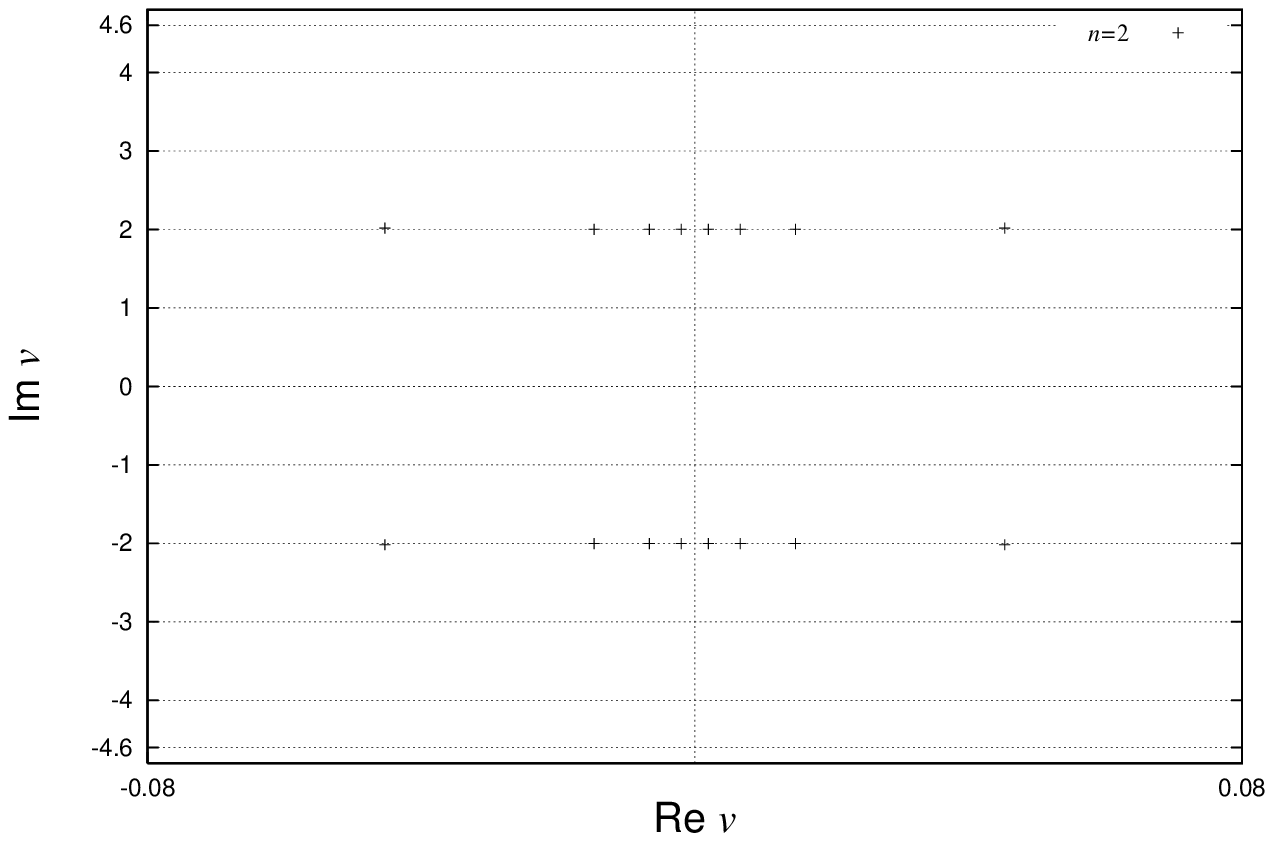,height=7cm}}
\caption{Location of zeros of $T^{(1)}_{n-1}(u,v)$ for
$n=2,3,4,5$, $p_0=\frac{24}{5}=4.8$, $u=-0.01$, $N=16$ (first),
$N=32$ (second) and $u=0.01$, $N=16$ (third). 
The zeros are located on an almost straight line 
$\Im v=\pm n\,\,\hbox{mod }2p_0$. The deviation from the line
 is $\sim 10^{-2}$.}
\label{zeros1}
\end{figure} 
\begin{figure}[p]
  \centerline{\epsfig{file=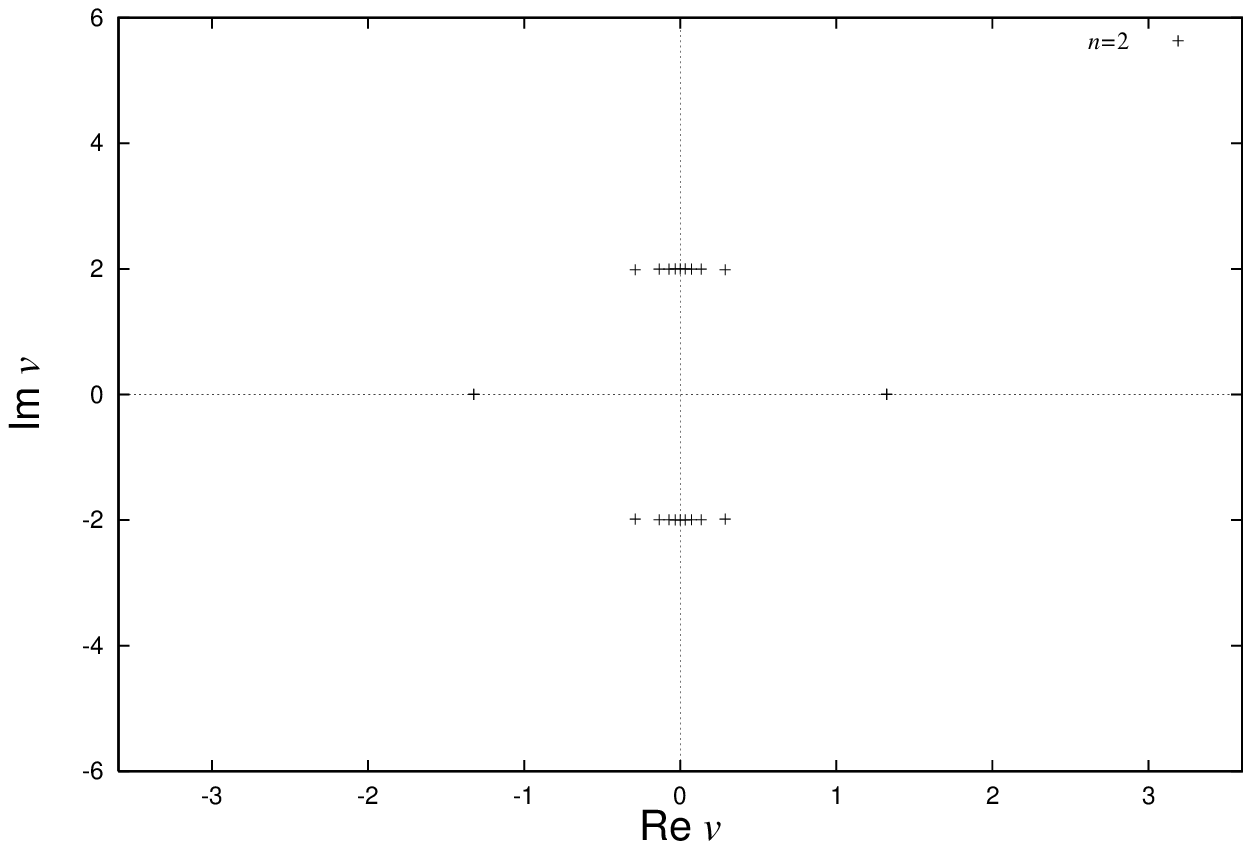}}
  \centerline{\epsfig{file=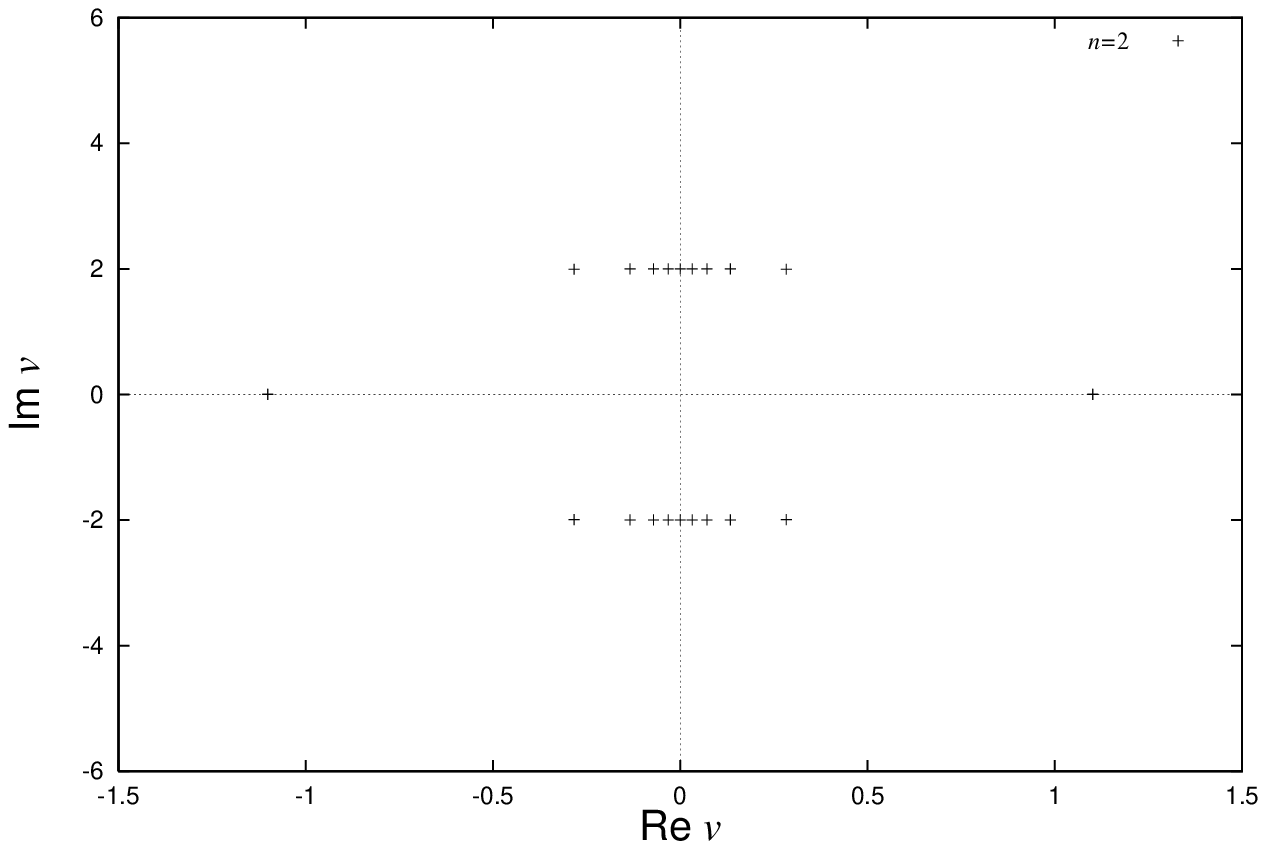}}
\caption{Location of zeros of $T^{(2)}_{n-1}(u,v)$ (upper)
and $T^{(3)}_{n-1}(u,v)$ (lower)
for $n=2,3,4,5$, $u=-0.1$, $p_0=5$, $N=20$.
$T^{(2)}_{n-1}(u,v)$ has two real zeros for $n\le 4$,
which are absent for $n=5$.
$T^{(3)}_{n-1}(u,v)$ has two real zeros for $n\le 4$ and
a double zero at $v=0$ for $n=5$. 
All the other zeros are located on an almost straight line 
$\Im v=\pm n\,\,\hbox{mod }2p_0$.}
\label{zeros23}
\end{figure} 
\begin{figure}[p]
  \centerline{\epsfig{file=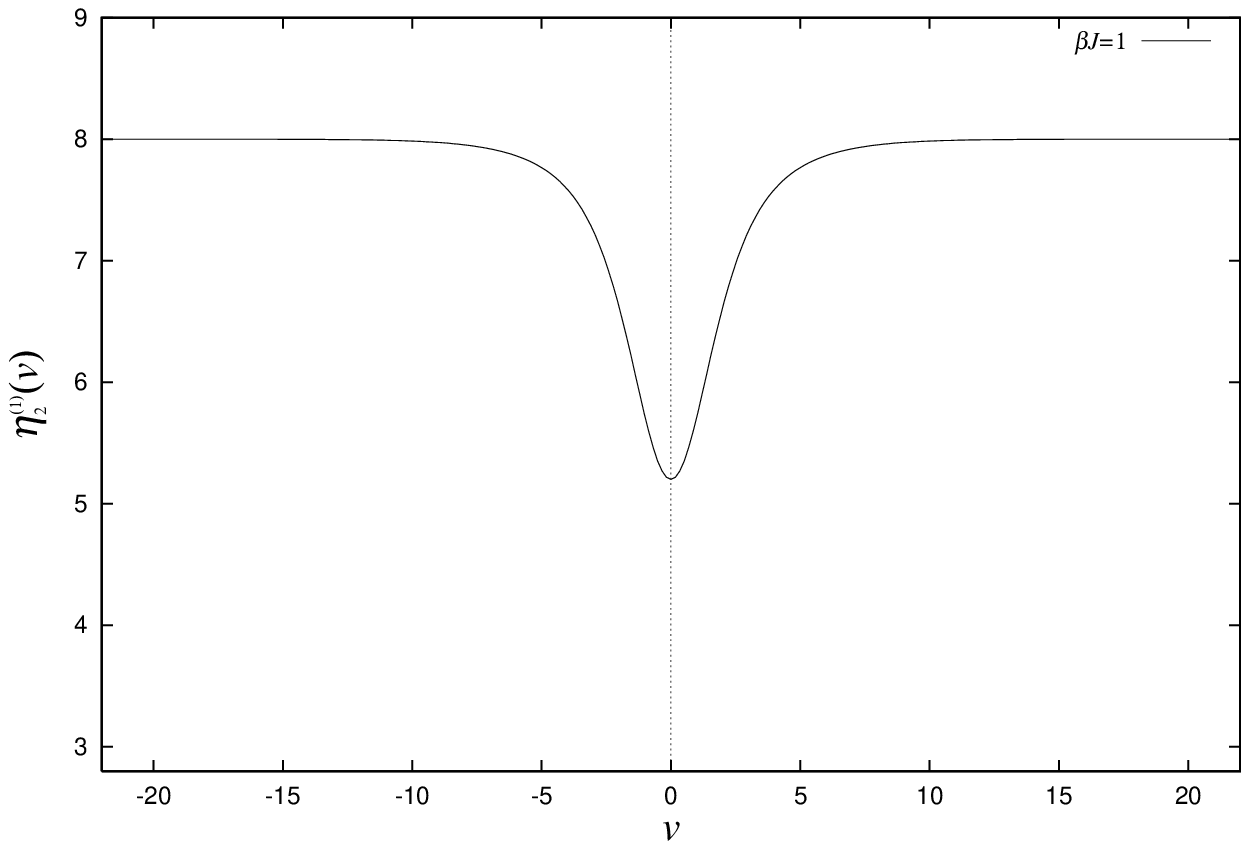,height=7.5cm}}
  \centerline{\epsfig{file=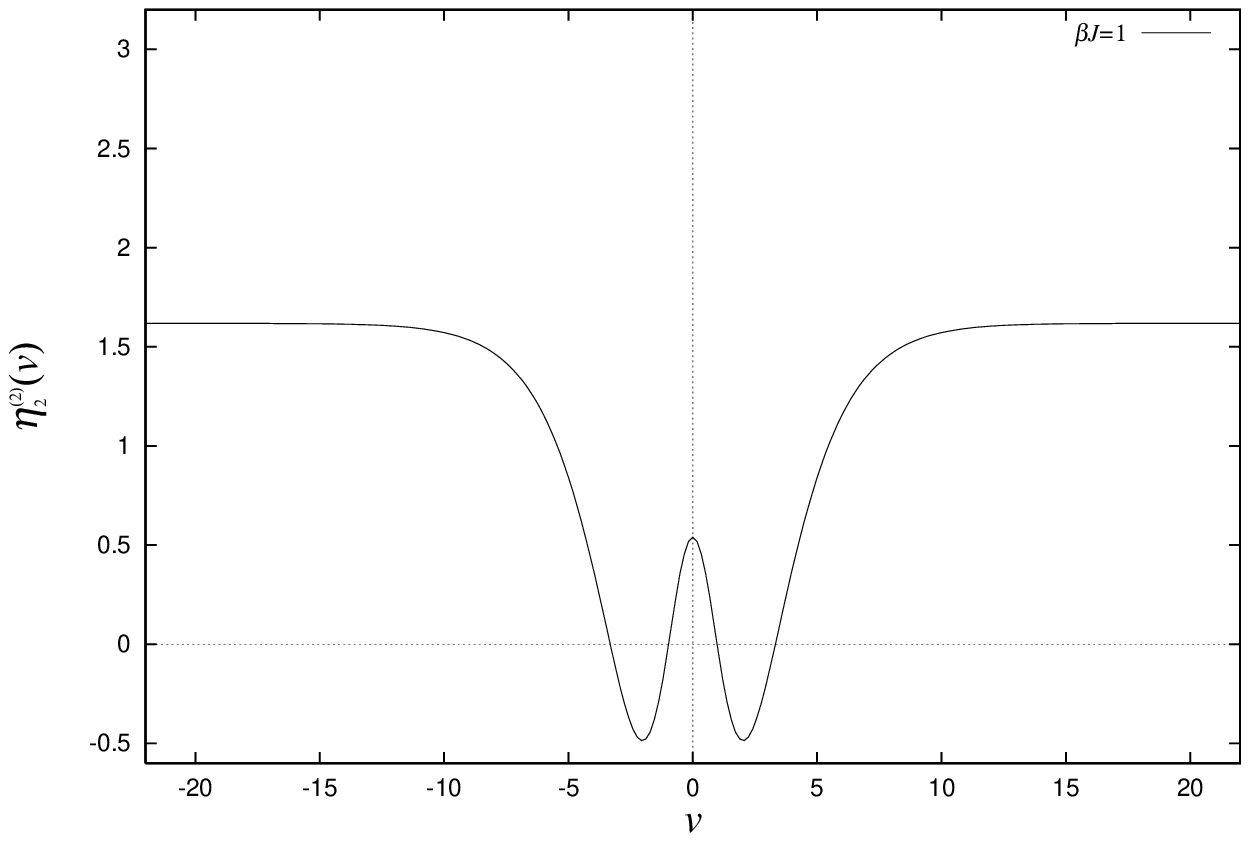,height=7.5cm}}
  \centerline{\epsfig{file=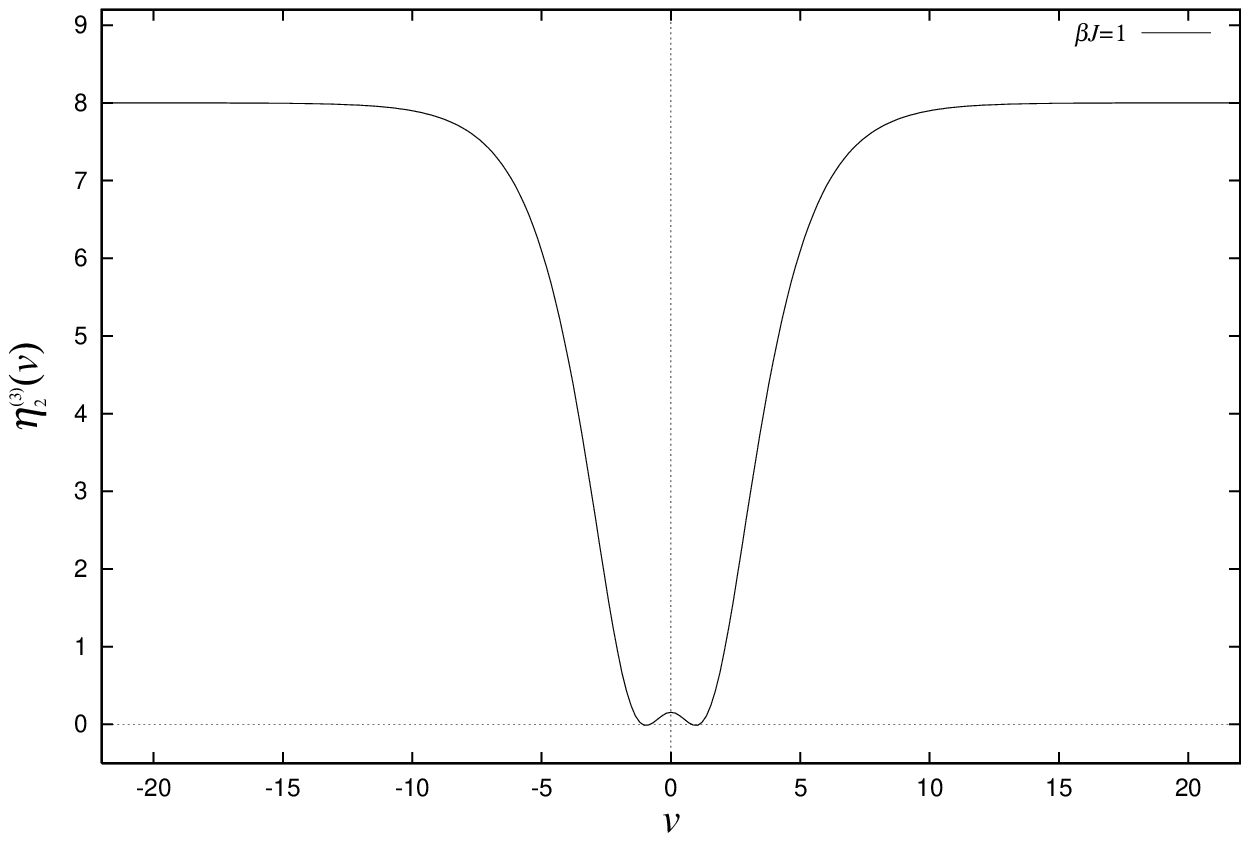,height=7.5cm}}
\caption{The function $\eta^{(k)}_2(v)$ for $v$ real, 
$p_0=5$ and $k=1,2,3$.}
\label{etag23}
\end{figure} 

\end{document}